\documentclass[12pt,a4paper]{book}

\usepackage{hyperref}
\usepackage{graphicx}
\usepackage{url}
\usepackage{amsthm}
\usepackage{newlfont}
\usepackage{setspace}
\usepackage{longtable}
\usepackage{ltxtable}
\usepackage{multirow}
\usepackage{amsmath}
\usepackage{amsfonts}
\usepackage{amssymb}

\usepackage{verbatim}

\usepackage{fancyhdr}

\input xypic
\xyoption{all}
\xyoption{poly}

\newtheorem{theo}{Theorem}

\newtheorem{cor}{Corollary}

\newtheorem{prop}{Proposition}
\newtheorem{fig}{Fig.}

\numberwithin{equation}{chapter} 
\numberwithin{theo}{section}
\numberwithin{prop}{section} 
\numberwithin{lem}{section}
\numberwithin{cor}{section} 
\numberwithin{rem}{section} 
\numberwithin{fig}{section}

\newcommand{\N}{\ensuremath{\mathbb N}}
\newcommand{\Z}{\ensuremath{\mathbb Z}}
\newcommand{\R}{\ensuremath{\mathbb R}}
\newcommand{\C}{\ensuremath{\mathbb C}}
\newcommand{\Q}{\ensuremath{\mathbb Q}}
\newcommand{\F}{\ensuremath{\mathbb F}}

\newcommand{\K}{\ensuremath{\mathbb K}}
\newcommand{\E}{\ensuremath{\mathbb E}}

\begin{document}
\pagenumbering{arabic}
\begin{titlepage}
\begin{singlespace}
\begin{center}
\begin{doublespace}
{\Huge{{\bf Paley Graphs and Their Generalizations}}}
\end{doublespace}
\vspace{0.5in}
By\\
    \vspace{0.5in}
  {\large{{\bf Ahmed Noubi Elsawy}}}\\
    \vspace{1.0in}
    A thesis submitted to\\
	Heinrich Heine University\\ 
	D\"{u}sseldorf, Germany\\
   for the Degree of Master of Science\\
\vspace{1.0in}
Supervisors\\
\vspace{0.5in}
{\large{{\bf Prof. Dr. F. Grunewald\quad\quad Prof. Dr. E. Klimenko}}}
\end{center}

\vfill
  \hfill
\parbox{2.5in}{Department of Mathematics\\
Heinrich Heine University\\
D\"{u}sseldorf\\
2009}
\end{singlespace}
\end{titlepage}

\thispagestyle{empty}
\begin{center}

  MASTER OF SCIENCE THESIS\\
  OF\\
  Ahmed Elsawy
\end{center}
\vfill
  APPROVED:\\\vspace{-1.0em}\\\\
\null\hfill\begin{tabular}{rl}
  THESIS COMMITTEE:\hspace{0.5in}\mbox{}\\\\
  MAJOR PROFESSOR & \rule[-0.1in]{3in}{0.4pt}\hspace{-3in}\\
  & \shortstack[l]{\rule{0in}{0.4in}\\\rule[-0.1in]{3in}{0.4pt}}\\
  & \shortstack[l]{\rule{0in}{0.4in}\\\rule[-0.1in]{3in}{0.4pt}}\\
  & \shortstack[l]{\rule{0in}{0.4in}\\\rule[-0.1in]{3in}{0.4pt}\\
                DEAN OF THE GRADUATE SCHOOL}\\
\end{tabular}
\vfill
\begin{center}
  HEINRICH HEINE UNIVERSITY\\
  D\"{U}SSELDORF\\
  2009  
\end{center}

\newpage
\pagestyle{empty}
\begin{center}
\Huge \bf Acknowledgment
\end{center}
\vspace{1cm}
I would like to offer my deepest gratitude and thankfulness to Prof. Dr. F. Grunewald for suggesting the
topics of this thesis, for his kind supervision, and for his invaluable help during the progress of present work.

I owe my deepest gratitude and special thankfulness to Prof. Dr. E. Klimenko who offered me a great
help in understanding the topics and to get this document into its present shape. Many thanks for her stimulating discussion, sincere advice, incisive comments, and careful reading of the manuscript.

Naturally, all the mistakes that can still be found in this thesis are mine. If you ever find a word
that is spelled correctly, it must have been Prof. Dr. E. Klimenko dropping me a line.

I want also to thank my family for their support and for keeping me going.\\
\vspace{0.5cm}

\hspace{\stretch{1}} A. Elsawy

\newpage
\thispagestyle{empty}
\vspace*{10pt}

\tableofcontents
\vspace{10pt}
\noindent \textbf{Bibliography}\hspace{\stretch{1}}\textbf{43}

\pagestyle{fancy} 
\fancyhead[OR]{{\it\leftmark}}
\fancyhead[OL]{{\it\thepage}}
\fancyhead[ER]{{\it\thepage}}
\fancyhead[EL]{{\it\rightmark}}
\fancyfoot[CO,CE]{}
\addtolength{\headheight}{3pt}

\chapter*{Introduction}
  \addcontentsline{toc}{chapter}{Introduction}
\thispagestyle{plain}

\indent Paley graphs are named after Raymond Paley (7 January 1907 -- 7 April 1933). He was born in Bournemouth, England. He won a Smith's Prize in 1930 and was elected a fellow of Trinity College, Cambridge, where he showed himself as one of the most brilliant students among a remarkable collection of fellow undergraduates. 
\\[10pt]
\indent In Paley graphs, finite fields form their sets of vertices. So to understand Paley graphs we will start our work with classification of finite fields and study their properties. We will show that any finite field $\F$ has $p^n$ elements, where $p$ is prime and $n \in \N$. Moreover, for every prime power $p^n$ there exists a field with $p^n$ elements, and this field is unique up to isomorphism.

In the last section of Chapter \ref{Chap1} we will see how one can construct the finite field for any prime power $p^n$, and we will give the explicit construction of the fields of 9, 16, and 25 elements.
\\[10pt]
\indent To construct a Paley graph, we fix a finite field and consider its elements as vertices of the Paley graph. Two vertices are connected by an edge if their difference is a square in the field. In the first section of Chapter \ref{Chap2} 
we will give some basic definitions and properties from graph theory which we will use in te study of the Paley graphs.
In the second section we will give the definition of the Paley graph and we will give examples of the Paley graphs of order 5, 9, and 13. Finally, we will study some important properties of the Paley graphs. 

In particular, we will show that the Paley graphs are connected, symmetric, and self-complementary.In \cite{W.P1}, Peisert proved that the Paley graphs of prime order are the only self-complementary symmetric graphs of prime order; furthermore, in \cite{W.P2}, he proved that any self-complementary and symmetric graph is isomorphic to a Paley graph, a $\mathcal{P}^*$-graph, or the exceptional graph $G(23^2)$ with $23^2$ vertices.

Also we will show that the Paley graph of order $q$ is $\frac{q-1}{2}$-regular, and every two adjacent vertices have $\frac{q-5}{4}$ common neighbors, and every two non-adjacent vertices have $\frac{q-1}{4}$ common neighbors, which means that the Paley graphs are strongly regular with parameters ($q,\ \frac{q-1}{2},\ \frac{q-5}{4},\ \frac{q-1}{4}$), see also \cite{Bollobas}.

\begin{displaymath}
\def\objectstyle{\scriptscriptstyle}
\xy /r6pc/:, {\xypolygon17"B"{\bullet}}, "B1"!{+U*++!L{\textrm{\normalsize0}}}
,"B2"!{+RD*+!LD{\textrm{\normalsize1}}},"B3"!{+RD*+!LD{\textrm{\normalsize2}}}, 
"B4"!{+LDD*++!D{\textrm{\normalsize3}}},"B5"!{+LD*++!D{\textrm{\normalsize 4}}},
"B6"!{+LDD*++!D{\textrm{\normalsize5}}},"B7"!{+LD*+!DR{\textrm{\normalsize6}}},
"B8"!{+LDD*+!DR{\textrm{\normalsize7}}},"B9"!{+LD*++!R{\textrm{\normalsize8}}},
"B10"!{+LDD*++!R{\textrm{\normalsize9}}},"B11"!{+LD*+!UR{\textrm{\normalsize10}}},
"B12"!{+LDD*+!UR{\textrm{\normalsize11}}},"B13"!{+LD*+!UR{\textrm{\normalsize12}}},
"B14"!{+LD*+!U{\textrm{\normalsize13}}},"B15"!{+LD*+!UL{\textrm{\normalsize14}}},
"B16"!{+LD*+!UL{\textrm{\normalsize15}}},"B17"!{+LD*+!L{\textrm{\normalsize16}}},
"B3";"B1"**@{-},"B3";"B5"**@{-},"B5";"B7"**@{-},"B7";"B9"**@{-},"B9";"B11"**@{-},
"B11";"B13"**@{-},"B13";"B15"**@{-},"B15";"B17"**@{-},"B17";"B2"**@{-},"B2";"B4"**@{-}
,"B4";"B6"**@{-},"B6";"B8"**@{-},"B8";"B10"**@{-},"B10";"B12"**@{-},"B12";"B14"**@{-},
"B14";"B16"**@{-},"B16";"B1"**@{-},
"B1";"B5"**@{-},"B9";"B5"**@{-},"B9";"B13"**@{-},"B13";"B17"**@{-},"B17";"B4"**@{-}
,"B4";"B8"**@{-},"B8";"B12"**@{-},"B12";"B16"**@{-},"B16";"B3"**@{-},"B3";"B7"**@{-}
,"B7";"B11"**@{-},"B11";"B15"**@{-},"B15";"B2"**@{-},"B2";"B6"**@{-},"B6";"B10"**@{-}
,"B10";"B14"**@{-},"B14";"B1"**@{-},
"B1";"B9"**@{-},"B9";"B17"**@{-},"B17";"B8"**@{-},"B8";"B16"**@{-},"B16";"B7"**@{-}
,"B7";"B15"**@{-},"B15";"B6"**@{-},"B6";"B14"**@{-},"B14";"B5"**@{-},"B5";"B13"**@{-}
,"B13";"B4"**@{-},"B4";"B12"**@{-},"B12";"B3"**@{-},"B3";"B11"**@{-},"B11";"B2"**@{-}
,"B2";"B10"**@{-},"B10";"B1"**@{-}
\endxy
\end{displaymath}

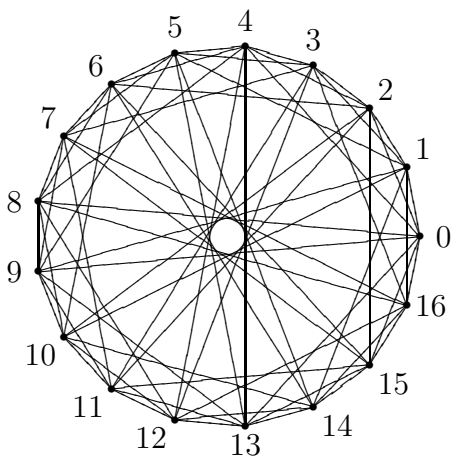
\begin{center}{\bf Fig. 0.1. }The Paley graph of order $17$
\end{center}

In \cite{ontheAdjpaley}, Ananchuen and Caccetta proved that for every 3-element subset $S$ of the vertices of the Paley graph with at least 29 vertices, and for every subset $T$ of $S$, there is a vertex $x \notin S$ which is joined to every vertex in $T$ and to no vertex in $S \setminus T$; that is, the Paley graphs are 3-existentially closed.
\\[10pt]
\indent Paley graphs are generalized by many mathematicians. In the first section of Chapter \ref{Chap3} we will see three examples of these generalizations and some of their basic properties. 

In \cite{ontheAdjGpaley}, Ananchuen introduced two of these generalizations. The cubic Paley graphs, in which pairs of elements of a finite field are connected by an edge if and only if they differ in a cubic residue, and the quadruple Paley graphs, in which pairs of elements of a finite field are connected by an edge if and only if they differ in a quadruple residue.  The third generalization is called the generalized Paley graphs, in this family of graphs, pairs of elements of a finite field are connected by an edge if and only if their difference belongs to a subgroup $S$ of the multiplicative group of the field. This generalization is given by Lim and Praeger in \cite{onGenPal}.

In the second section of Chapter \ref{Chap3} we will define a new generalization of the Paley graphs, in which pairs of elements of a finite field are connected by an edge if and only if there difference belongs to the $m$-th power of the multiplicative group of the field, for any odd integer $m > 1$, and we call them the $m$-Paley graphs. 

Since the cubic Paley graphs are 3-Paley graphs, we can say that the cubic Paley graphs are a special case of the family of $m$-Paley graphs. Also, we will give some examples of this family.

In the third section we will show that the $m$-Paley graph of order $q$ is complete if and only if $\gcd (m, q-1)= 1$ and when $d= \gcd (m, q-1) > 1$, the $m$-Paley graph is $\frac{q-1}{d}$-regular. 

Also we will prove that the $m$-Paley graphs are symmetric but not self-complementary. In particular, $m$-Paley graphs are not in the Peisert's list. Since strongly regular graphs must be self-complementary, we see that the $m$-Paley graphs are not strongly regular.

We will show also that the $m$-Paley graphs of prime order are connected but the $m$-Paley graphs of order $p^n, \ n > 1 $ are not necessary connected, for example they are disconnected if $\gcd (m,p^n-1)= \frac{p^n-1}{2}$.

\chapter{Finite Fields}
\label{Chap1} 
This chapter provides an introduction to some basic properties of finite fields and their structure. This introduction will be useful for understanding the properties of Paley graphs in which the elements of a finite field represent the set of vertices. 
\begin{section}{Basic definitions and properties}
\label{Sec1}

\noindent {\bf Definition:}
A {\it field} \ $\F$ is a set  of at least two elements, with two operations $\oplus$  and  $\ast$, for which the following axioms are satisfied:

\begin{enumerate}
  \item The set \ $\F$ under the operation $\oplus$  forms an abelian group (whose identity is denoted by $0$).
  \item The set \ $\F^{*} = \F \setminus \lbrace0\rbrace$ under the operation $\ast$ forms an abelian group (whose identity is denoted by 1).
  \item Distributive law: For all $a$, $b$, $c \in \F$, we have $(a \oplus b)\ast c = (a \ast c) \oplus (b \ast c)$.
\end{enumerate}
\noindent Note that it is important to allow $0$ to be an exceptional element with no inverse, because if $0$ had an inverse then $1 = 0^{-1} \ast 0 = 0$, it would follow that $x = x\ast 1 = 0 $ for all $ x \in \F $; therefore, $\F$ would consist of only one element 0. 
\\[10pt]
{\bf\large{Examples}}
\begin{enumerate}
\item $\Q$, $\R$, and $\C$ are fields with respect to the usual addition and multiplication.
\item The subring  $\Q[i] := \lbrace a + bi \in \C$ : $a$, $b \in \Q \rbrace$  of  $\C$  is a field, called the field of Gaussian rationals.
\item The ring of integers modulo $p$, $\Z_{p}$, is a field if $p$ is a prime number.
\item A commutative division ring (a ring in which every nonzero element has a multiplicative inverse) is a field.
\item A finite integral domain (a commutative ring with no zero divisors) is a field.

Indeed, let $R$ be an integral domain and $0 \neq a \in R$. The map $x \rightarrow ax$, $x \in R$, is injective because $R$ is an integral domain $(a(x_{1} - x_{2}) = 0 \Leftrightarrow (x_{1} - x_{2}) = 0)$. If $R$ is finite, the map is surjective as well, so that $ax = 1$ for some $x$, i.e., 
every nonzero element $a$ has an inverse.
\item  The quotient ring $R/M$ such that $R$ is a commutative ring and $M$ is a maximal ideal, is a field.

Indeed, since $R$ is commutative, the ring $R/M$ is commutative, where $(x +M)(y +M) = xy +M = yx +M = (y+M)(x +M)$. It also has an identity $1_{R/M} = 1_{R} + M$. Moreover, $1_{R} \notin M$ because if $1_{R} \in M$, then $M = R$, which contradicts the definition of maximal ideal. In order to show that $R/M$ is a field, it remains to prove that every nonzero element $x+M \in R/M$ has an inverse. So we fix $x \notin M$ and consider the set 
$I = M +xR =\{ m+xr : m \in M, \ r \in R \} $. First we need to check that $I$ is an ideal. Let $m_{1},m_{2} \in M$ and $r_{1}, r_{2}, s \in R$,
then $0 = 0 + x0 \in I$, $s((m_1 + xr_1) - (m_2 + xr_2)) = s(m_1-m_2) + s(x(r_1-r_2)) = s(m_1-m_2) + xs(r_1-r_2) \in I$. 
Hence $I$ is an ideal. Now $M \subsetneq I$ because $x \in I$ but $x \notin M$. Since $M$ is maximal, it follows that $I = R$, 
and in particular $1_R \in I$. So there exist $m \in M$ and $y \in R$ such that $1_R = m+xy$. 
Then $(x+M)(y+M) = xy+M =(1_R - m) +M = 1_R +M$, so $x+M$ has an inverse in $R/M$.
\end{enumerate}
\noindent {\bf Definition:} The {\it characteristic} of a field  $\F$ (denoted by Char $\F$) is the smallest positive integer $n$ such that $n1=0$, where $n1$ is an abbreviation for $1 + 1 + \cdots +1 (n \textrm{ ones})$. If $n1$ is never $0$, we say that \ $\F$ has characteristic $0$.
\\[10pt]
\noindent Note that the characteristic of a field can never be equal to $1$, since $1= 1 \ast 1 \neq 0$. If Char $\F \neq 0$, then Char $\F$ must be a prime number. For if Char $\F = n = rs$ where $r$ and $s$ are positive integers greater than $1$, then $(r1)(s1) = n1 = 0$, so either $r1$ or $s1$ is $0$, which contradicts the minimality of $n$.
\\[10pt]
\noindent {\bf Definition:}
If \ $\F$ and $\E$ are fields and \ $\F \subseteq \E$, we say that $\E$ is an {\it extension} of \ $\F$ and \ $\F$ is a {\it subfield} of \ $\E$, and we write \ $\F \leq \E$.
\\[10pt]
\noindent Note that if \ $\F \leq \E$ then we can consider $\E$ as a vector space over $\F$, because if we consider the elements of $\E$ as vectors and the elements of $\F$ as scalars, then the axioms of a vector space are satisfied. The dimension of this vector space is called the degree of the extension and is denoted by $[\E:\F]$. If $[\E:\F]= n$ with $n < \infty$, we say that $\E$ is a {\it finite extension} of $\F$.
\\[10pt]
\noindent {\bf Definition:}
A minimal subfield $\F_{p}$ of a field $\F$ with Char $\F = p$ is called a {\it prime field}.
\\[10pt]
\noindent Note that the only subfield of the prime field $\F_{p}$ is $\F_{p}$ itself. Let Char $\F= p$ then $\lbrace1,1+1,1+1+1,\dots,0=1+1+\cdots+1(p \textrm{ ones} )\rbrace$ form a subfield of $\F$ which is a prime field isomorphic to $\Z_{p}$, i.e., every finite field with characteristic $p$ has a prime subfield which is isomorphic to $\Z_{p}$.
\\[10pt]
\noindent {\bf Definition:}
A {\it nonzero polynomial }$f(x)$ of degree $m$ over a field $\F$ is an expression of the form
\begin{center}
$f(x) = f_{0} + f_{1} x + f_{2} x^{2} + \cdots + f_{m} x^{m}$
\end{center}
where $f_{i} \in \F \textrm{, }0 \leq i \leq m \textrm{, and }f_{m} \neq 0$. The degree of $f(x)$ is denoted by $\textrm{deg}{\hskip 0.7mm} f(x)$. The polynomial $f(x) = 0 $ is called the {\it zero polynomial} and the set of all polynomials over $\F$ is denoted by $\F[x]$.
\\[10pt]
\noindent {\bf Definition:}
A nonzero polynomial \ $f(x)= f_{0} + f_{1} x + f_{2} x^{2} + \cdots + f_{m} x^{m}$ over \ $\F$  with $f_{m}=1$ is called {\it monic}.
\\[10pt]
\noindent Note that the set $\F[x]$ is a ring, its additive identity is the zero polynomial $f(x) = 0$ and its multiplicative identity is $f(x) = 1$. However $\F[x]$ is not a field, because the polynomials of degree greater than 0 have no inverse.
\\[10pt]
\noindent {\bf Definition:}
A nonzero polynomial  $f(x) \in \F[x]$ is called {\it irreducible} over $\F$ if $\textrm{deg}{\hskip 0.7mm}f(x) \geq 1$ and $f(x) = g(x)h(x)$ with $g(x)$, $h(x) \in \F[x]$ gives either $\textrm{deg}{\hskip 0.7mm} g(x) = 0$ or $\textrm{deg}{\hskip 0.7mm}h(x) = 0$.
\\[10pt]
\noindent In other words, up to a nonzero constant factor the only divisors of $f(x)$ are $f(x)$ itself and $1$.
\\[10pt]
\noindent {\bf Definition:}
Let  $f(x)$ be a nonzero polynomial of $\F[x]$. A finite extension  $\E$ of  $\F$ is called the {\it splitting field} of  $f(x)$ if   $\E$ is the smallest extension field of  $\F$ in which  $f(x)$ can be written as
$$\lambda (x - \alpha_{1}) \cdots (x - \alpha_{k})$$ 
for some $\alpha_{1} , \dots , \alpha_{k} \in \E$ and $\lambda$ in $\F$.
\\[10pt]
Note that the splitting field $\E$ can be written as $\E = \F[\alpha_{1} , \dots , \alpha_{k}]$ which denotes the field generated by $\alpha_{1} , \dots , \alpha_{k}$ over $\F$, because $\F[\alpha_{1} , \dots , \alpha_{k}]$ is the smallest field containing $\F$ and $(\alpha_{1} , \dots , \alpha_{k})$.
\\[10pt]
\noindent {\it Example: }The field of complex numbers is the splitting field of $x^{2} + 1$ over the field of real numbers.
\end{section}

\begin{section}{Classification of the finite fields}
\label{sec2}

\begin{theo}
  The number of elements of a finite field  $\F$ is equal to $p^{n}$, where $p$ is prime and $n \in \N$.
\end{theo}

\proof
Let $\textrm{Char}{\hskip 0.5mm}\F = p$ where $p$ is a prime number, then $\F_{p}$ is a subfield of $\F$, so we can consider $\F$ as a vector space over $\F_{p}$. Since $\F$ is finite, the dimension of $\F$ over $\F_{p}$ is finite. 
\\[10pt]
\indent Let $[\F:\F_{p}] =n$ then there exists a basis {$v_{1}, v_{2}, \dots , v_{n}$} in $\F$. With respect to this basis every element $x \in \F$ can be written uniquely as
$$x = a_{1} v_{1} + a_{2} v_{2} + \cdots + a_{n} v_{n}, \textrm{ where } a_{1}, a_{2},\dots ,a_{n} \in \F_{p}.$$
\indent Since $\vert \F_{p} \vert= p$, every $a_{i}$, $i = 1,2, \dots ,n$, can be equal to one of $p$ values. Therefore, $\F$ has $p^{n}$ elements. \hspace{\stretch{1}} $\Box$ 
\\[10pt]
\noindent Now we will consider the following, more difficult, question: Does a field with $p^{n}$ elements exist, for each prime $p$ and each $n \in \N$?  This question will be answered affirmatively in Theorem \ref{p^n}. In order to prove this theorem, we need Theorems \ref{theo EX&UN} and \ref{theo f&f'}.

\begin{theo}
\label{theo EX&UN}
Let  $f(x)$ be an irreducible polynomial with degree $\geqslant 1$ in  $\F[x]$. Then a splitting field for $f(x)$ over  $\F$ exists and any two such splitting fields are isomorphic.
\end{theo}

\proof
First we will prove that there exists an extension field of  $\F$ in which $f(x)$ has a root. 
\\[10pt]
\indent Since $f(x)$ is irreducible, the principal ideal $(f(x))$ is a maximal ideal in the ring $\F[x]$. 
To prove this, let $I$ be an ideal in $\F[x]$ with $(f(x)) \subsetneq I \subseteq \F[x]$ and let $g(x) \in I \setminus (f(x))$. Since $f(x)$ is irreducible and $f(x) \nmid g(x)$, we have $(f(x),g(x))= 1$. Therefore, there exist two polynomials $h(x)$, $k(x) \in \F[x]$ with $f(x)h(x)+g(x)k(x) = 1$. We see that $ 1 \in I$ and thus  $I = \F[x]$, so $(f(x))$ is a maximal ideal in the commutative ring $\F[x]$. We conclude that $\F[x]/(f(x))$ is a field. 
\\[10pt]
\indent Consider the homomorphism
$$\sigma: \F \rightarrow \F[x]/(f(x)) \textrm{ defined by } \sigma(a) = \overline{a} = a + (f(x)),\textrm{ }a\in \F,$$
\noindent which is injective. Indeed, $\F$ is a field, which means that $\textrm{Ker} \sigma$ equals $\F$ or $(0)$  but $\textrm{Ker} \sigma \neq \F$ because $\overline{1} \neq \overline{0}$. So $\textrm{Ker} \sigma = (0) $ and $\F$ is isomorphic to $\sigma (\F)$. 
\\[10pt]
\indent Now we can consider $\F[x]/(f(x))$ as a finite extension of $\F$, which has a root $\alpha = x+(f(x))$ of $f(x)$. Thus for a field $\F$ and an irreducible polynomial $f(x)$ over $\F$ there exists an extension in which $f(x)$ has a root.
\\[10pt]
\indent Now we will prove by induction on $\textrm{deg} {\hskip 0.7mm} f(x)$ the existence of the splitting field of $f(x)$. 

If $\textrm{deg} {\hskip 0.7mm} f(x) = 1$ then $f(x)$ has one root. Thus, as we have shown, $\F$ has an extension field $\E$ which contains a root $\alpha$ of $f(x)$, so the field $\F[\alpha]$ is the splitting field for $f(x)$. 

Now assume that for each irreducible polynomial $g(x)$ of degree $m < n$ there exists a splitting field of $g(x)$. 

Let $\textrm{deg} {\hskip 0.7mm} f(x) = n$, we have shown that there exists an extension field $\E_{1}$ of $\F$ with root $\alpha_{1}$ of $f(x)$. So in $\E_{1}$, $f(x)$ can be written as 
\begin{center} 
$f(x)=(x-\alpha_{1})g(x)$, with $\textrm{deg} {\hskip 0.7mm} g(x) = n-1$.
\end{center}

By induction hypothesis, there exists a splitting field in which $g(x)$ can be written as 
\begin{center}
$g(x)=c(x-\alpha_{2}) \cdots (x-\alpha_{n})$, with $c \in \F$. 
\end{center}

Thus the field $\E = \F[\alpha_{1},\alpha_{2},\dots,\alpha_{n}]$ is the splitting field of $f(x)$ over $\F$ in which all the roots of $f(x)$ are contained. 
\\[10pt]
\indent We have proved the existence of the splitting field of $f(x)$, let us now prove the uniqueness up to isomorphism. 

Let $\E$ and $\K$ be two splitting fields of $f(x)$, then there exists a nontrivial isomorphism $\theta: \F \longrightarrow \F$ such that $\E = \F[\alpha_{1},\dots,\alpha_{n}]$ and $\K = \theta\F[\beta_{1},\dots,\beta_{n}]$ with $f(x) = c_{1}(x-\alpha_{1}) \cdots (x-\alpha_{n})$ over $\E$ and $\theta f(x) = c_{2}(x-\beta_{1}) \cdots (x-\beta_{n})$ over $\K$ and $c_{1} , c_{2} \in \F$. 

We will prove by induction on $\textrm{deg} {\hskip 0.7mm} f(x)$, that $\theta$ can be extended to an isomorphism $\overline{\theta}: \E \longrightarrow \K$. 

If $\textrm{deg} {\hskip 0.7mm} f(x)= 1$ then $\E = \F[\alpha]$ and $\K = \theta\F[\beta]$. Thus the mapping $\theta_{1}: \E \longrightarrow \K$ defined by $\theta_{1}(\alpha)= \beta$ and $\theta_{1}(c)= \theta(c)$ for all $c \in \F$, is an isomorphism. 

If $\textrm{deg} {\hskip 0.7mm} f(x)= n$, then for a root $\alpha_{1}$ of $f(x)$ and a root $\beta_{1}$ of $\theta f(x)$, there exists an isomorphism  $\theta_{1}:\F[\alpha_{1}] \longrightarrow \theta\F[\beta_{1}]$. So we can write $f(x)$ as 

\noindent $f(x)= c_{1}(x-\alpha_{1})g(x)$ over $\F[\alpha_{1}]$ and $\theta f(x)= c_{2}(x-\beta_{1})\theta_{1}g(x)$ over $\theta\F[\beta_{1}]$. 

Since $\textrm{deg} {\hskip 0.7mm} g(x)= n-1 =$ $\textrm{deg} {\hskip 0.7mm} \theta_{1}g(x)$, by induction hypothesis, there exists an isomorphism $\theta_{2} : \F[\alpha_{2},\dots,\alpha_{n}] \longrightarrow \theta\F[\beta_{2},\dots,\beta_{n}]$ such that, after a permutation of $\beta_{2},\dots,\beta_{n}$, $\theta_{2}(\alpha_{i})=\beta_{i}$ for $i= 2,\dots,n$. 
\\[10pt]
Thus we can write $g(x)$ and $\theta g(x)$ as 

$g(x)= c_{3}(x-\alpha_{2})\cdots(x-\alpha_{n})$ over $\F[\alpha_{2},\dots,\alpha_{n}]$, 

$\theta g(x)= c_{4}(x-\beta_{2})\cdots(x-\beta_{n})$ over $\theta\F[\beta_{2},\dots,\beta_{n}]$. 
\\[10pt]
So $f(x)= c_{1}(x-\alpha_{1})g(x)= c(x-\alpha_{1})\cdots(x-\alpha_{n})$ over $\F[\alpha_{1},\dots,\alpha_{n}]$ and $\theta f(x)= c_{2}(x-\beta_{1})\theta_{1}g(x)=\theta(c)(x-\beta_{1})\cdots(x-\beta_{n})$ over $\theta\F[\beta_{1}, \dots ,\beta_{n}]$ such that $c=c_{1}c_{3}$ and $\theta(c)= c_{2}c_{4}$. Then $\theta$ can be extended to $\overline{\theta}:\E \longrightarrow \K$ such that $\overline{\theta}(\alpha_{i})=\beta_{i}$ for $i= 1,\dots,n$ and $\overline{\theta}(c) =\theta(c)$ for all $c \in \F$.  \hspace{\stretch{1}} $\Box$

\begin{theo}
\label{theo f&f'}
Let $\F$ be a field and $f(x)$ be a polynomial over $\F$. Then $f(x)$ has no roots with multiplicity $\geq 2$ if and only if the greatest common divisor of a polynomial $f(x)$ and $f^{'}(x)$ has degree $0$. i.e., $f(x)$ and $f^{'}(x)$ has no common root.
\end{theo}

\proof
Assume that the greatest common divisor of $f(x)$ and $f^{'}(x)$ has degree $0$, we will prove that $f(x)$ has no roots with multiplicity $\geq 2$ by contradiction. Let $f(x)$ have at least one root with multiplicity $r \geq 2$, then we can write $f(x)$ as 
$$f(x) = c(x-a_{1})^{r}(x-a_{2})\cdots(x-a_{n-r})$$
where $n = \textrm{deg} {\hskip 0.7mm} f(x) , c \in \F$ and $a_{1},a_{2},\dots,a_{n-r}$ are all roots of $f(x)$ such that $a_{1}$ is repeated $r$ times. Then
\begin{center}
$f^{'}(x) = cr(x-a_{1})^{r-1}(x-a_{2}) \cdots (x-a_{n-r})+c(x-a_{1})^{r}(x-a_{3}) \cdots (x-a_{n-r})+ \cdots +c(x-a_{1})^{r}(x-a_{2}) \cdots (x-a_{n-r-1}) = c(x-a_{1})^{r-1}[r(x-a_{2}) \cdots (x-a_{n-r})+(x-a_{1})(x-a_{3}) \cdots (x-a_{n-r})+ \cdots +(x-a_{1})(x-a_{2}) \cdots (x-a_{n-r-1})]$.
\end{center} 
It follows that $(x-a_{1})^{r-1}\mid f(x)$ and $(x-a_{1})^{r-1}\mid f^{'}(x)$, so the greatest common divisor of $f(x)$ and $f^{'}(x)$ has degree $\geq r-1 \neq 0$ because $r \geq 2 $.
\\[10pt]
\indent Now assume that $f(x)$ has no roots with multiplicity $\geq 2$, we prove that the greatest common divisor of a polynomial $f(x)$ and $f^{'}(x)$ has degree $0$. We can write $f(x)$ as
$$f(x) = c(x-a_{1})^{r}(x-a_{2})\cdots(x-a_{n})$$
where $n = \textrm{deg} {\hskip 0.7mm} f(x), c \in \F$ and $a_{1},a_{2},\dots,a_{n}$ are all roots of $f(x)$.
Then
$f^{'}(x)= c(x-a_{2})\cdots(x-a_{n})+c(x-a_{1})(x-a_{3})\cdots(x-a_{n})+\cdots+c(x-a_{1})(x-a_{2})\cdots(x-a_{n-1})=c[(x-a_{2})\cdots(x-a_{n})+(x-a_{1})(x-a_{3})\cdots(x-a_{n})+\cdots+(x-a_{1})(x-a_{2})\cdots(x-a_{n-1})]$. 

It follows that $c$ is the only common divisor of $f(x)$ and $f^{'}(x)$. \hspace{\stretch{1}} $\Box$

\begin{theo}
\label{p^n}
 For every prime $p$ and $n \in \N$ there is a field with $p^{n}$ elements.
\end{theo}
\proof 
Consider the polynomial $f(x)= x^{p^{n}}-x$ over a field $\F$. If $f(x)$ is irreducible, then by Theorem \ref{theo EX&UN} there exists a splitting field $\E$ of $f(x)$ which is unique up to isomorphism. If $f(x)$ is not irreducible, set $\E = \F$. Let $K$ be the set of all zeros of $f(x)$ in $\E$, then $K = \{ a \in \E : a^{p^{n}}=a \}$. 
\\[10pt]
\indent Since $f^{'}(x) = p^{n}x^{(p^{n}-1)}-1$, $f(x)$ and $f^{'}(x)$ has no common zero, where $f(a)=0 \Rightarrow a^{p^{n}}=a \Rightarrow a^{p^{n}-1}=1 \Rightarrow f^{'}(a)=p^{n}a^{p^{n}-1}-1=p^{n}-1 \neq 0$. 

It follows from Theorem \ref{theo f&f'} that all zeros of $f(x)$ are distinct. Thus the set $K$ has $p^{n}$ elements.
\\[10pt]
\indent We will prove that $K = \E$. In order to show that, we need to prove that $K$ is a field, in which each element is a root of $f(x)$, then $K$ is a splitting field of $f(x)$ with $p^{n}$ elements. 

It is clear that $K \subseteq \E$ and by using Theorem \ref{theo EX&UN} $K$ isomorphic to $\E$. Thus $K = \E$.
\\[10pt]
\indent Clearly $0,1 \in K$. Let $a , b \in K$, to prove that $K$ is closed under addition, we need to prove that $(a+b)^{p^{n}}= a^{p^{n}}+b^{p^{n}}= a+b$. In the polynomial expansion 
\begin{displaymath}
(a+b)^{p^{n}} =  \sum_{i=0}^{p^{n}} \binom{p^{n}}{i} a^{i} b^{p^{n}-i} : \binom {p^{n}}{i}= \dfrac{p^{n}!}{i! (p^{n}-i)!}	\end{displaymath} 
we can see that all binomial coefficients are divisible by $p^{n}$ except the first and the last, and since the finite field $\E$ has  $\textrm{Char} {\hskip 0.7mm} \E = p$, all binomial coefficients are $0$ except the first and the last. That is, $$(a+b)^{p^{n}}= a^{p^{n}}+b^{p^{n}}=a+b.$$ 

Clearly, $K$ is closed under multiplication or $(ab)^{p^{n}}= a^{p^{n}}b^{p^{n}}$, because $\E$ is commutative. Since $\forall \ a \neq 0 ,\ a^{p^{n}} = a$, we have $(a^{-1})^{p^{n}} = (a^{p^{n})^{-1}} = a^{-1} $. So the inverse of any element in $K$ belongs to $K$. 

Thus $K$ is subfield of $\E$ or $K = \E$ is a field with $p^{n}$ element. \hspace{\stretch{1}} $\Box$
\\[10pt]
\noindent {\bf Notation} : We will denote by $\F_{p^n}$ the field with $p^n$ elements. 
\\[10pt]
Clearly $(\F_{p^n})^*=\F_{p^n} \setminus \{0\}$ and $(\F^{*}_{p^n})^2=\{\ a^2 : a \in \F^{*}_{p^n} \}$ form a group under multiplication. In general, $(\F^{*}_{p^n})^m=\{\ a^m : a \in \F^{*}_{p^n}$, $m \in \N\}$, form a group under multiplication. 
\\[10pt]
\begin{theo}
\label{TH:cyclic}
The multiplicative group $\ \F^*_{p^n}$ is cyclic.
\end{theo}
\proof
Let $p^n \geqslant 3$ and $h=p^n-1=p^{r_1}_1p^{r_2}_2\cdots p^{r_m}_m$ be the prime factorization of $|\F^*_{p^n}|$. Let $a_i$, for every $1 \leqslant i \leqslant m$, be an element in $\F_{p^n}$ with $a^{\frac{h}{p_i}}_i\neq 1$. To prove the existence of $a_i$, consider the polynomial $x^{\frac{h}{p_i}}-1$ which has at most $\frac{h}{p_i}$ roots in $\F_{p^n}$. Since $\frac{h}{p_i} < h$, it follows that there exists $a_i \in \F_{p^n}$ with $a^{\frac{h}{p_i}}_i-1 \neq 0$ or $a^{\frac{h}{p_i}}_i \neq 1$. 

Now set $b_i=a^{{h}/{p^{r_i}_i}}_i$, it follows that $b_i \neq 1$ and $b^{p^{r_i}}_i = 1$, then the order of $b_i$ must be in the form $p^{s_i}_i$ with $1 \leq s_i \leq r_i$.  Since $$b^{p^{r_i-1}_i}_i = (a^{{h}/{p^{r_i}_i}}_i)^{p^{r_i-1}_i}= a^{\frac{h}{p_i}}_i \neq 1,$$ the order of $b_i$ is $p^{r_i}_i$. 

Now we will prove that $b= b_1b_2\cdots b_m$ is a generator of the group $\ \F^*_{p^n}$, that is the order of $b$ is $h$. We know that $b^h=b^h_1b^h_2\cdots b^h_m=1$. 

Assume that the order of $b$ is not $h$, then the order of $b$ is a proper divisor of $h$. So the order of $b$ is a divisor of at least one of the $m$ integers $\frac{h}{p_i}$ and $1 \leq i \leq m$, say $\frac{h}{p_1}$. Thus 
$$b^{h/p_1}=b^{h/p_1}_1b^{h/p_1}_2\cdots b^{h/p_1}_m=1.$$ Since $p^{r_i}_i | \frac{h}{p_1}$ for every $2 \leq i \leq m$, it follows that $b^{h/p_1}_i=1$ for every $2 \leq i \leq m$. Thus $b^{h/p_1}_1=1$, which is impossible, because the order of $b_1$ is $p^{r_1}_1 \nmid \frac{h}{p_1}$. 

So the order of $b$ is $h$, in other words, the multiplication group $\ \F^*_{p^n}$ is cyclic.
\hspace{\stretch{1}} $\Box$
\\[10pt]
\noindent We remark here that $\F_{p^n}$ with $n > 1$ is never $\Z_{p^n}$, in the following section we construct some of such fields.
\end{section}

\begin{section}{Construction of the finite fields}
\label{sec3}

To construct a field with $p^{n}$ elements, we use an irreducible monic polynomial $f(x) \in \Z_{p}[x]$ with $\textrm{deg} {\hskip 0.7mm} f(x) = n$. The elements of the field $\Z_{p}[x]/(f(x))$ can be written in the form
$$ a_{0} + a_{1} x + a_{2} x^{2} + \cdots + a_{n-1} x^{n-1}, \textrm { where }  a_{i} \in \Z_{p} \textrm { for all } \ i = 0,1,\dots,n-1.$$

Since there are $p$ possible values for each $a_{i}$, the field $\Z_{p}[x]/(f(x))$ has $p^{n}$ elements. 

To find an irreducible polynomial, we list all possible monic polynomials of degree $n$ ($p^{n}$ possible monic polynomials), which is not always an easy process especially for large $n ,p$. 

Clearly, any polynomial without a constant term is not irreducible ($x$ is a factor), so these $p^{n-1}$ polynomial will not be considered. 

For each of the remaining $p^{n} - p^{n-1}$ polynomials, we could substitute one by one all the field elements for $x$. If none of these substitutions is equal to zero, the polynomial is irreducible (i.e., it has no root in the field). 
\\[10pt]
\indent Let $a$ be a zero of the chosen polynomial, then the elements of $\Z_{p}[x]/(f(x))$ can be written in its vector form representation using the basis $$\{1, a, a^{2},\dots, a^{n-1}\}.$$ 

\indent We can also generate a multiplicative representation of the field by using the fact that the multiplicative group of the field is cyclic. So if we can find a primitive element (i.e., a generator of the cyclic group), we will have a representation of the elements. 
\\[10pt]
{\it Example 1}: We will construct a field of \ $16 = 2^{4}$ \ elements, here $p=2 , n=4$. We start with a field of order $2$ which is \ $\Z_{2}= \{0,1\}$ and an irreducible polynomial over $\Z_{2}$ of degree $4$. We can easily list all possible polynomials of degree $4$ over $\Z_{2}$. There are $16$ of them :
\begin{flushleft}
$x^4 ,\quad x^4 + 1,\quad x^4 + x,\quad x^4 + x^2,\quad x^4 + x^3,\quad x^4 + x + 1,\quad x^4 + x^2 + 1$,\\
$x^4 + x^3 + 1,\quad x^4 + x^2 + x,\quad x^4 + x^3 + x,\quad x^4 + x^3 + x^2,\quad x^4 + x^2 + x + 1$,\\
$x^4 + x^3 + x + 1,\quad x^4 + x^3 + x^2 + 1,\quad x^4 + x^3 + x^2 + x,\quad x^4 + x^3 + x^2 + x +1$.
\end{flushleft}
Every polynomial without a constant term has root $0$. So we will consider just the following polynomials
\begin{flushleft}
$x^4 + 1,\quad x^4 + x + 1,\quad x^4 + x^2 + 1,\quad x^4 + x^3 + 1,\quad x^4 + x^2 + x + 1$,\\ 
$x^4 + x^3 + x + 1,\quad x^4 + x^3 + x^2 + 1,\quad x^4 + x^3 + x^2 + x +1$.
\end{flushleft}
In this set, every polynomial with even number of terms has root $1$. So we will consider just the following polynomials
\begin{flushleft}
$x^4 + x + 1,\quad x^4 + x^2 + 1,\quad x^4 + x^3 + 1,\quad x^4 + x^3 + x^2 + x +1$.
\end{flushleft}
All of them have no roots in $\Z_{2}$, then all of them are irreducible polynomials over $\Z_{2}$. 

Consider one of them, say $x^4 + x + 1$. Let $a$ be a root of $x^4 + x + 1$, then the elements of the field $\Z_{2}[x]/(x^4 + x + 1)$ can be obtained by two methods: 

The first method is additive, in which we construct all linear combinations of $1, a, a^2$ and $a^3$. They are :
\begin{flushleft}
$0,\quad 1,\quad a,\quad a^2,\quad a^3,\quad a + 1,\quad a^2 + 1,\quad a^3 + 1,\quad a^2 + a,\quad a^3 + a$,\\ 
$a^3 + a^2,\quad a^2 + a + 1,\quad a^3 + a + 1,\quad a^3 + a^2 + 1,\quad a^3 + a^2 + a$,\\ 
$a^3 + a^2 + a + 1$.
\end{flushleft}
The second method is multiplicative. Since $a^4 = -a - 1 = a + 1$, we can write down the powers of $a$ as the following: 
\begin{displaymath} 
\begin{array}{lllllllll}
a^1 &=& a & a^2 &=& a^2 & a^3 &=& a^3 \\
a^4 &=& a + 1&a^5 &=& a^2 + a&a^6 &=& a^3 + a^2\\
a^7 &=& a^3 + a + 1&a^8 &=& a^2 + 1&a^9 &=& a^3 + a\\
a^{10} &=& a^2 + a + 1& a^{11} &=& a^3 + a^2 + a&a^{12} &=& a^3 + a^2 + a + 1\\
a^{13} &=& a^3 + a^2 + 1& a^{14} &=& a^3 + 1& a^{15} &=& 1,
\end{array}
\end{displaymath}
which means that $a$ is a generator of the cyclic group 
$$(\Z_{2}[x]/(x^4 + x + 1))^{*} = \Z_{2}[x]/(x^4 + x + 1) \setminus \{0\}.$$ 

Notice also that the terms on the right are all the possible terms that can be written as linear combinations of the basis $\{1, a, a^2, a^3\}$ over $\Z_{2}$. When working with finite fields it is convenient to have both of the above representations, since the terms on the left are easy to multiply and the terms on the right are easy to add. 
\\[10pt]
\indent Now suppose we had chosen a root of the second irreducible polynomial $x^4 + x^2 + 1$, say, $b$. We would then have  $b^4= b^2 + 1$ and the powers of $b$ will be 
\begin{displaymath}
\begin{array}{lllllllll}
b^1 &=& b&b^2&=& b^2& b^3 &=& b^3\\
b^4 &=& b^2 + 1& b^5 &=& b^3 + b& b^6 &=& b^4 + b^2 = b^2 + 1 + b^2 = 1,
\end{array}
\end{displaymath}
which means that $b$ cannot be a generator of the group $(\Z_{2}[x]/(x^4 + x^2 + 1))^{*}$.
\\[10pt]
{\it Example 2}: Now we will construct a field of \ $9 = 3^{2}$ \ elements, that is, $p=3 , n=2$. We start with a field of order $3$ which is \ $\Z_{3}= \{0, 1, 2\}$ and an irreducible polynomial over $\Z_{3}$ of degree $2$. We can easily list all possible monic polynomials over $\Z_{3}$. They are :
\begin{flushleft}
$x^2,\quad x^2 + 1,\quad x^2 + 2,\quad x^2 + x,\quad x^2 + 2x,\quad x^2 + x + 1,\quad x^2 + x + 2,\quad x^2 + 2x + 1,$\\
$x^2 + 2x + 2$.
\end{flushleft}
Every polynomial without a constant term has root $0$. So we will consider just the following polynomials
\begin{flushleft}
$x^2 + 1
,\quad x^2 + 2 
,\quad x^2 + x + 1,\quad 
x^2 + x + 2
,\quad x^2 + 2x + 1,\quad 
x^2 + 2x + 2.$ 
\end{flushleft}
In this set, $x^2 + 2$ and $x^2 + x + 1$ have root $1$ and $x^2 + 2x + 1$ has root $2$. So we will consider just the following polynomials
\begin{center}
$x^2 + 1,\quad x^2 + x + 2,\quad x^2 + 2x + 2$.
\end{center}
All of them have no roots in $\Z_{3}$, then all of them are irreducible polynomials over $\Z_{3}$. 

Consider one of them, say $x^2 + x + 2$. Let $a$ be a root of $x^2 + x + 2$, then the elements of the field $\Z_{2}[x]/(x^2 + x + 2)$ can be obtained by two methods: 

The first method is additive, in which we construct all linear combinations of $1$ and $a$. They are :
\begin{center}
$0,\quad 1,\quad 2,\quad a,\quad 2a,\quad a + 1,\quad a + 2,\quad 2a + 1,\quad 2a + 2$.
\end{center}

The second method is multiplicative. Since $a^2 = -a - 2 = 2a + 1$, we can write out the powers of $a$ as follows:
\begin{displaymath}
\begin{array}{llllllllllll}
a^1 &=& a
&a^2 &=& 2a + 1
&a^3 &=& 2a + 2
&a^4 &=& 2\\
a^5 &=& 2a&
a^6 &=& a + 2
&a^7 &=& a + 1
&a^8 &=& 1.
\end{array}
\end{displaymath}

In other words, $a$ is a generator of the cyclic group 
$$(\Z_{3}[x]/(x^2 + x + 2))^{*} = (\Z_{2}[x]/(x^2 + x + 2) )\setminus \{0\}.$$

\noindent {\it Example 3}: Finally, we will construct a field of \ $25 = 5^{2}$ \ elements. Here, $p=5 , n=2$, a field of order $5$ is \ $\Z_{5}= \{0, 1, 2, 3, 4\}$. To find an irreducible polynomial over $\Z_{5}$ of degree $2$, we list all possible monic polynomials of degree $2$ over $\Z_{5}$:
\begin{flushleft}
$x^2,\quad x^2 + 1,\quad x^2 + 2,\quad x^2 + 3,\quad x^2 + 4,\quad x^2 + x,\quad x^2 + 2x,\quad x^2 + 3x$,\\
$x^2 + 4x,\quad x^2 + x + 1,\quad x^2 + x + 2,\quad x^2 + x + 3 ,\quad x^2 + x + 4,\quad x^2 + 2x + 1$,\\ 
$x^2 + 2x + 2,\quad x^2 + 2x + 3,\quad x^2 + 2x + 4,\quad x^2 + 3x + 1,\quad x^2 + 3x + 2,\quad x^2 + 3x + 3$,\\
$x^2 + 3x + 4,\quad x^2 + 4x + 1,\quad x^2 + 4x + 2,\quad x^2 + 4x + 3,\quad x^2 + 4x + 4$.
\end{flushleft}
Every polynomial without a constant term has root $0$,\\
$x^2 + 4$, $x^2 + x + 3$, $x^2 + 2x + 2$ and $x^2 + 3x + 1$ have root $1$,\\
$x^2 + 1$ ,$x^2 + 3x$, $x^2 + x + 4$, $x^2 + 2x + 2$ and $x^2 + 4x + 3$ have root $2$,\\
$x^2 + 1$, $x^2 + 2x$, $x^2 + x + 3$, $x^2 + 3x + 2$ and $x^2 + 4x + 4$ have root $3$,\\
and $x^2 + 4$, $x^2 + x$, $x^2 + 2x + 1$, $x^2 + 3x + 2$ and $x^2 + 4x + 3$ have root $4$.\\
So we will consider the remaining polynomials:
\begin{flushleft}
$x^2 + 2,\quad x^2 + 3,\quad x^2 + x + 1,\quad x^2 + x + 2,\quad x^2 + 2x + 3,\quad x^2 + 2x + 4$,\\ 
$x^2 + 3x + 3,\quad x^2 + 3x + 4,\quad x^2 + 4x + 1,\quad x^2 + 4x + 2$.
\end{flushleft}
All of them have no roots in $\Z_{5}$, then all of them are irreducible polynomials over $\Z_{5}$. 

Consider one of them, say $x^2 + 2$. Let $a$ be a root of $x^2 + 2$, then the elements of the field $\Z_{5}[x]/(x^2 + 2)$ are all linear combinations of $1$ and $a$. They are :
\begin{flushleft}
$0,\quad 1,\quad2,\quad3,\quad 4,\quad a,\quad2a,\quad3a,\quad4a,\quad a + 1,\quad a + 2,\quad  a + 3,\quad  a + 4$, \\
$ 2a + 1,\quad 2a + 2,\quad 2a + 3,\quad 2a + 4,\quad3a + 1,\quad 3a + 2, \quad3a + 3,\quad 3a + 4$,\\
$4a + 1,\quad 4a + 2,\quad4a + 3,\quad4a + 4$.
\end{flushleft}
\end{section}
\chapter{Paley Graph}
\label{Chap2}
In this chapter we will give some basic definitions and properties of graph theory and we will study in details Paley graphs and some of its properties.
\begin{section}{Basic definitions and properties}
\noindent {\bf Definition:}
A {\it graph} $G$ is a pair $( V, E )$ of sets satisfying $E \subseteq P_2(V)$, where $P_2(V)$ is the set of all subsets of $V$ with two elements. The elements of\ $V$ are called {\it vertices} and the elements of $E$ are called {\it edges}.
\\[10pt]
\noindent Note that the set of vertices of a graph $H = ( W, F)$ is denoted by $V(H)$ and the set of edges is denoted by $E(H)$. An edge $e = \{x, y\}$ is sometimes written as $xy$. 
\\[10pt]
\noindent {\bf Definition:}
The {\it order} of a graph $G$ is the number of its vertices and is denoted by $|G|$. A graph $G$ is called a {\it finite graph} or an {\it infinite graph} depending on the order of $G$. If $|G| = 0$ then $G$ is called the {\it empty graph}.
\\[10pt]
\noindent In order to draw a graph $G$, we can represent its vertex set $V(G)$ by dots, we join two of these dots by a line if and only if the two corresponding vertices form an edge in $E(G)$.
\begin{displaymath}
\def\objectstyle{\scriptscriptstyle}
\xy /r2.8pc/:, {\xypolygon6"A"{~>{}{\bullet}}}, "A1"!{+U*+!L{\textrm{\normalsize1}}}
,"A2"!{+RD*++!D{\textrm{\normalsize2}}},"A6"!{+RD*+!UL{\textrm{\normalsize6}}}
,"A3"!{+RD*++!DR{\textrm{\normalsize3}}}, "A4"!{+LDD*+!R{\textrm{\normalsize4}}},
"A5"!{+LD*+!UR{\textrm{\normalsize5}}}, "A4";"A2"**@{-}, "A1";"A3"**@{-}, "A5";"A2"**@{-}
, "A1";"A2"**@{-}, "A1";"A5"**@{-}, "A2";"A3"**@{-}, "A1";"A4"**@{-}
\endxy
\end{displaymath}
\begin{fig}
\begin{center}
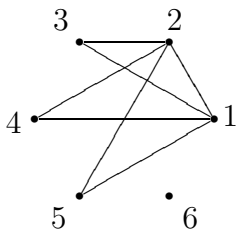
A graph $G$ with $V(G)=\{1,2,3,4,5,6\}$ and $E(G)=\{\{1,2\},\{1,3\},\{1,4\},\{1,5\},\{2,3\},\{2,4\},\{2,5\}\}$.
\end{center}
\end{fig}
\noindent {\bf Definition:}
Two vertices $x$, $y$ of $G$ are {\it adjacent}, or {\it neighbors}, if $xy$ is an edge of $G$, and the set of all neighbors of $x$ is denoted by $N(x)$. 

The order of the set $N(x)$ is called the {\it degree} of $x$ and is denoted by $d(x)$. We say that a vertex $x$ is {\it isolated} if  $d(x) = 0$.
\\[10pt]
\noindent Note that for each graph $G = ( V, E )$ we have \begin{displaymath}
 2|E(G)| = \sum_{x\in V(G)} d(x).\end{displaymath}

\noindent {\bf Definition:}
A graph $G$ is called {\it k-regular} if all the vertices of $G$ have the same degree $k$. 

If $|V(G)|= n$ and $G$ is an ($n$-1)-regular graph then $G$ is called a {\it complete} graph. We will denote by $K_n$ the complete graph on $n$ vertices.
\\[10pt]
\noindent Note that in a $k$-regular graph $G$ we have $2|E(G)| = k |V(G)|$, it follows that $k$ or $|V(G)|$ is even.
\\[10pt]
\noindent {\bf Definition:}
A {\it path} of length $n$ in a graph $G$ is the sequence 
$$x_{1}e_{1}x_{2} \cdots e_{n-1}x_{n} \textrm{ with } x_{i} \in V(G),\ e_{i}=\{x_{i},x_{i+1}\}\in E(G)$$ for all $i \in \{ 1, 2, \dots, n-1\}$, and $x_{i}\neq x_{j}$ for all $i\neq j$.

Let $x_{1}e_{1}x_{2}e_{2} \cdots x_{n-1}e_{n-1}x_{n}$ be a path then the sequence 
$$x_{1}e_{1}x_{2}e_{2} \cdots x_{n-1}e_{n-1}x_{n}e x_{1} \textrm{ with } e=\{x_{n},x_{1}\}\in E(G)$$ 

\noindent is called a {\it cycle} of length $n$ which will be denoted by $C_n$.
\\[10pt]
\noindent As an example, in figure \ref{fig2} we see the cycle and the complete graph of $5$ vertices. 

\begin{displaymath}
\def\objectstyle{\scriptscriptstyle}
\xy /r3pc/:, {\xypolygon5"B"{\bullet}}, "B1"!{+U*+!L{\textrm{\normalsize1}}}
,"B2"!{+RD*++!D{\textrm{\normalsize2}}}
,"B3"!{+RD*+!R{\textrm{\normalsize3}}}, "B4"!{+LDD*+!UR{\textrm{\normalsize4}}},
"B5"!{+LD*+!UL{\textrm{\normalsize5}}}
\endxy
\quad\quad\quad
\xy /r3pc/:, {\xypolygon5"A"{\bullet}}, "A1"!{+U*+!L{\textrm{\normalsize1}}}
,"A2"!{+RD*++!D{\textrm{\normalsize2}}} ,"A3"!{+RD*+!R{\textrm{\normalsize3}}}, 
"A4"!{+LDD*+!UR{\textrm{\normalsize4}}}, "A5"!{+LD*+!UL{\textrm{\normalsize5}}}, 
"A4";"A1"**@{-}, "A1";"A3"**@{-}, "A5";"A2"**@{-}, "A2";"A4"**@{-}, "A3";"A5"**@{-}
\endxy
\end{displaymath}
\begin{center}
$C_5$\quad\quad\quad\quad\quad
\quad\quad\quad\quad $K_5$\end{center}
\begin{fig}\label{fig2}
\begin{center}
\end{center}
\end{fig}
\noindent {\bf Definition:}
Let $G$ be a $k$-regular graph with $|G|=n$. If there are two integers $\lambda, \mu$ such that 

every two adjacent vertices have $\lambda$ common neighbors and 

every two non-adjacent vertices have $\mu$ common neighbors, 

\noindent then $G$ is called a {\it {\bf s}trongly {\bf r}egular {\bf g}raph} with parameters $(n, k, \lambda, \mu)$ and is denoted by {\bf srg}$(n, k, \lambda, \mu)$.
\\[10pt]
\indent Clearly, every strongly regular graph is regular, but not vice versa. For example $C_6$ is 2-regular but not a strongly regular graph.
\\[10pt]
\indent Note that in any strongly regular graph srg($n,k,\lambda,\mu$), its parameters are related by 
\begin{equation}\label{eq1} \mu (n-k-1)=k(k-\lambda-1). \tag{*} \end{equation}
\indent In order to show that, consider a vertex $x \in V(G)$. We remind that $N(x)$ is the set of all neighbors of $x$. Let $N^{'}(x)$ be the set of all non-adjacent vertices of $x$, then $n = |N(x)|+1+|N^{'}(x)|$. We will prove that both sides of (\ref{eq1}) are equal to the number of edges between $N(x)$ and $N^{'}(x)$. 
\\[10pt]
\indent Since $G$ is $k$-regular, $|N^{'}(x)|= n -|N(x)|-1 = n-k-1$. Let $y \in N^{'}(x)$, then $y$ and $x$ have $\mu$ common neighbors, that is $\mu$ equals the number of edges between $y$ and $N(x)$. It follows that $\mu (n-k-1)$ is the number of edges between $N(x)$ and $N^{'}(x)$. 
\\[10pt]
\indent Let $z \in N(x)$ then $z$ and $x$ have $\lambda$ common neighbors. Since $|N(x)|= k$, it follows that the number of neighbors of $z$ which are not adjacent to $x$ is equal to $k-\lambda -1$. Thus $k-\lambda -1$ is the number of edges between $z$ and $N^{'}(x)$, so $k(k-\lambda -1)$ equals the number of edges between $N(x)$ and $N^{'}(x)$. Therefore, (\ref{eq1}) is proved. 
\\[10pt]
\noindent {\bf Definition:}
A graph $G$ is called {\it connected } if every two vertices are connected by a path.
\\[10pt]
\indent Note that every complete graph is connected, regular, and strongly regular. Both strongly regular and regular graphs are not necessary connected and also connected graphs are not necessary complete, strongly regular, or regular. 

For example, a cycle is a connected graph which is not complete, a path is a connected graph which is not regular, and the graph $G$ with $V(G)=\{1,2,3,4\}$ and $E(G)=\{\{1,2\},\{3,4\}\}$ is strongly regular which is not connected. 

\begin{displaymath}
\def\objectstyle{\scriptscriptstyle}
\xy /r3pc/:, {\xypolygon4"B"{\bullet}}, "B1"!{+U*+!DL{\textrm{\normalsize1}}}
,"B2"!{+RD*+!DR{\textrm{\normalsize2}}}
,"B3"!{+RD*+!UR{\textrm{\normalsize3}}}, "B4"!{+LDD*+!UL{\textrm{\normalsize4}}}
\endxy
\quad\quad\quad
\xy /r3pc/:, {\xypolygon4"B"{~>{}\bullet}}, "B1"!{+U*+!DL{\textrm{\normalsize1}}}
,"B2"!{+RD*+!DR{\textrm{\normalsize2}}}
,"B3"!{+RD*+!UR{\textrm{\normalsize3}}}, "B4"!{+LDD*+!UL{\textrm{\normalsize4}}},
"B1"\PATH~={**@{-}}'"B2"'"B3"'"B4"
\endxy
\quad\quad\quad
\xy /r3pc/:, {\xypolygon4"B"{~>{}\bullet}}, "B1"!{+U*+!DL{\textrm{\normalsize1}}}
,"B2"!{+RD*+!DR{\textrm{\normalsize2}}}
,"B3"!{+RD*+!UR{\textrm{\normalsize3}}}, "B4"!{+LDD*+!UL{\textrm{\normalsize4}}}
,"B1";"B2"**@{-}, "B3";"B4"**@{-}
\endxy
\end{displaymath}
\begin{fig}
\begin{center}
\end{center}
\end{fig}
\noindent {\bf Definition:}
Let $G=(V, E) \textrm{ and } G^{'}=(V^{'}, E^{'})$ be two graphs. $G$ is {\it isomorphic} to $G^{'}$ if there is a bijection $f : V\longrightarrow V^{'}$ such that $xy \in E$ if and only if $f(x)f(y)\in E^{'}$; we denote this by $G \cong G^{'}$. 

An isomorphism from a graph $G$ to itself is called an {\it automorphism}. The set of all automorphismus of a graph $G$ form a group under composition, and it is denoted by $Aut(G)$. 
\\[10pt]
\noindent {\bf Definition:}
Let $G = (V, E)$ be a finite graph. The {\it complementary} graph of $G$ is a graph $\overline{G}$ with $V(\overline{G})=V(G)$ and $E(\overline{G})= P_2(V) \setminus E(G)$. That is, $xy \in E(\overline{G})$ if and only if $xy \notin E(G)$. 

A graph $G$ is called {\it self-complementary} if it is isomorphic to its complement.
\\[10pt]
For example $C_5$ is self-complementary, see figure 2.1.4. 

\begin{displaymath}
\def\objectstyle{\scriptscriptstyle}
\xy /r3pc/:, {\xypolygon5"B"{\bullet}}, "B1"!{+U*+!L{\textrm{\normalsize1}}}
,"B2"!{+RD*++!D{\textrm{\normalsize2}}}
,"B3"!{+RD*+!R{\textrm{\normalsize3}}}, "B4"!{+LDD*+!UR{\textrm{\normalsize4}}},
"B5"!{+LD*+!UL{\textrm{\normalsize5}}},
"B4";"B1"**@{.}, "B1";"B3"**@{.}, "B5";"B2"**@{.}, "B2";"B4"**@{.}, "B3";"B5"**@{.}
\endxy
\textrm{ \quad\quad\quad }
\xy /r3pc/:, {\xypolygon5"A"{~>{.}\bullet}}, "A1"!{+U*+!L{\textrm{\normalsize1}}}
,"A2"!{+RD*++!D{\textrm{\normalsize2}}} ,"A3"!{+RD*+!R{\textrm{\normalsize3}}}, 
"A4"!{+LDD*+!UR{\textrm{\normalsize4}}}, "A5"!{+LD*+!UL{\textrm{\normalsize5}}}, 
"A4";"A1"**@{-}, "A1";"A3"**@{-}, "A5";"A2"**@{-}, "A2";"A4"**@{-}, "A3";"A5"**@{-}
\endxy
\end{displaymath}
\begin{center}
$C_5\quad\quad\quad\quad\cong\quad\quad\quad\quad\overline{C_5}$
\end{center}
\begin{fig}
\begin{center}
\end{center}
\end{fig}
\noindent {\bf Definition:}
Let $G$ be a group and let $X$ be a non-empty set. We say that $G$ acts on $X$ if there is a map $ \phi:G\times X \rightarrow X$ such that the following conditions hold for all $x \in X$:

\begin{enumerate}
 \item $\phi(e,x)=x$ where $e$ is the identity element of $G$.
 \item $\phi(g,\phi(h,x))=\phi(gh,x) \ \forall g,h \in G$.
\end{enumerate}

In this case, $G$ is called a transformation group, $X$ is a called a $G$-set, and $\phi$ is called the group action.

The group action is called {\it transitive} (we also say that $G$ acts {\it transitively} on $X$) 
if for every $x,y \in X$, there exists $g \in G$ such that $\phi(g,x)= y$.
\\[10pt]
\noindent {\bf Definition:}
A graph $G$ is called {\it symmetric} if its automorphism group acts transitively on the vertices and edges.
\\[10pt]
\indent For example every cycle graph or complete graph is symmetric graph and every symmetric graph is a regular graph, but not vice versa.
\begin{displaymath}
\def\objectstyle{\scriptscriptstyle}
\xy /r3.5pc/:, {\xypolygon6"A"{\bullet}}, "A1"!{+U*+!L{\textrm{\normalsize1}}}
,"A2"!{+RD*+!D{\textrm{\normalsize2}}} ,"A3"!{+RD*+!D{\textrm{\normalsize3}}}, 
"A4"!{+LDD*+!R{\textrm{\normalsize4}}},"A5"!{+LDD*+!U{\textrm{\normalsize5}}},
"A6"!{+LDD*+!U{\textrm{\normalsize6}}}, 
"A1";"A4"**@{-},"A2";"A5"**@{-},"A3";"A6"**@{-}
\endxy
\textrm{ \quad\quad\quad }
\xy /r3.5pc/:, {\xypolygon7"A"{~>{}\bullet}}, "A1"!{+U*+!L{\textrm{\normalsize1}}}
,"A2"!{+RD*+!D{\textrm{\normalsize2}}} ,"A3"!{+RD*+!R{\textrm{\normalsize3}}}, 
"A4"!{+LDD*+!R{\textrm{\normalsize4}}},"A5"!{+LDD*+!U{\textrm{\normalsize5}}},
"A6"!{+LDD*+!U{\textrm{\normalsize6}}},"A7"!{+LDD*+!L{\textrm{\normalsize7}}}, 
"A1";"A2"**@{-},"A2";"A3"**@{-},"A3";"A1"**@{-},"A4";"A5"**@{-},"A5";"A6"**@{-},
"A6";"A7"**@{-},"A7";"A4"**@{-}
\endxy
\end{displaymath} 
\begin{fig}
\begin{center}
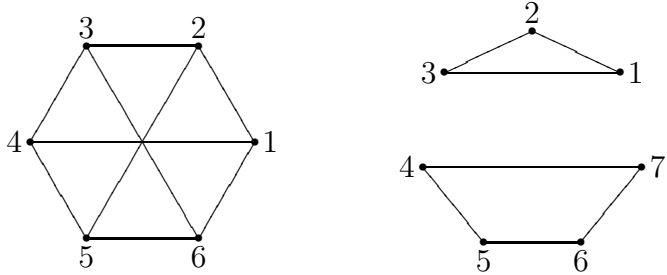
Symmetric graph which is regular and regular graph which is not symmetric.
\end{center}
\end{fig}

\end{section}
\begin{section}{Paley graphs}
\begin{subsection}{Definition and examples}
Before we define the Paley graph, we need the following definition.
\\[10pt]
{\bf Definition:} Let $q$ and $r$ be two positive integers with  gcd $(q,r)=1$, then $r$ is a {\it quadratic residue} of $q$ if and only if $x^2 \equiv r \ (\textrm{mod}\ q$) has a solution, and $r$ is a {\it quadratic nonresidue} of $q$ if and only if $x^2 \equiv r \ (\textrm{mod}\ q$) has no solution.
\\[10pt]
{\bf Definition:} Let $p$ be a prime number and $n$ be a positive integer such that $p^n \equiv 1\ (\textrm{mod}\ 4)$. The graph $P =(V, E)$ with $$V(P)=\F_{p^n} {\textrm \ and \ } E(P)=\{\{x,y\} : x,y \in \F_{p^n},\ x-y\in (\F^{*}_{p^n})^2\}$$ is called the {\it Paley graph} of order $p^n$.
\\[10pt]
\indent Note that the set $E(P)$ in the definition of Paley graph is well defined because $ x-y\in (\F^*_{p^n})^2$ if and only if $y-x\in (\F^*_{p^n})^2$. Since $x-y= -1(y-x)$, we need only to show that $-1 \in (\F^*_{p^n})^2$. 

We have $p^n \equiv 1\ (\textrm{mod}\ 4)$, so $4 \mid (p^n -1)$. Let $g$ be a generator of the group $\F^*_{p^n}$ then $p^n-1$ is the least positive integer such that $g^{p^n-1} = 1$. We can rewrite this as $g^{p^n-1}-1=(g^{\frac{p^n-1}{2}}-1)(g^{\frac{p^n-1}{2}}+1)=0$. Since $g^{\frac{p^n-1}{2}}$ cannot be equal to $1$, it follows that $g^{\frac{p^n-1}{2}}=(g^{\frac{p^n-1}{4}})^2 =-1$ which means that $g^{\frac{p^n-1}{4}}$ is a square root of $-1$.
\\[10pt]
\indent Note that if the Paley graph has prime order $p$, then we can consider the field of integers modulo $p$, $\Z_p$, as its vertex set. 
\\[10pt]
\indent However, we cannot consider $\Z_{p^n}$ with $n > 1$ as a vertex set of the Paley graph of order $p^n$, because as we have seen in the previous chapter, there exists a unique field $\F_{p^n}$ of order $p^n$ which is not $\Z_{p^n}$, such a field will represent the set of vertices of the Paley graph. To get this field $\F_{p^n}$ we need the construction that was developed in the previous chapter. 
\\[10pt]
\indent The list of integers which can be considered as an order of the Paley graph starts with $5$, $9$, $13$, $17$, $25$, $29$, $37$, $41$. In the following examples, we show the Paley graphs explicitly for the first three cases. 
\\[10pt]
{\it Example 1:}
The Paley graph of order $5$ is the cycle $C_5$. 
\\[10pt]
\indent In order to see that, let $P=(V, E)$ be the Paley graph of order $5$ then $V(P)=\Z_5=\{0,1,2,3,4\}$ and $(\Z^*_5)^2=\{1,4\}$, it follows that $$E(P)=\{\{0,1\},\{1,2\},\{2,3\},\{3,4\},\{4,0\}\}.$$
{\it Example 2:}
Let $P=(V, E)$ be the Paley graph of order $9=3^2$. Here we have $p=3$, $n=2$, then $V(P)=\F_{3^2}$, the field of order $9$, can be written as $$\F_{3^2}=\{0,1,2,a,2a,1+a,1+2a,2+a,2+2a\}\cong \Z_3[x]/(x^2+1)$$ where $a$ is a root of $x^2+1$. Since 
$$
\begin{array}{cllcllcll}
1^2&=&1,
&2^2&=&1,
&a^2&=&-1=2,\\
(2a)^2&=&2,&
(1+a)^2&=&2a,&
(1+2a)^2&=&a,\\
(2+a)^2&=&a,
&(2+2a)^2&=&2a,
\end{array}
$$
we have $(\F^*_{3^2})^2=\{1,2,a,2a\}$. Thus $E(G)=\{\{0,1\},\{0,2\},\{0,a\},\{0,2a\},$
$\{1,2\},\{1,1+a\},\{1,1+2a\},\{2,2+a\},\{2,2+2a\},\{a,1+a\},\{a,2+a\},\{a,2a\}$,
$\{2a,1+2a\},\{2a,2+2a\},\{1+a,2+a\},\{1+a,1+2a\},\{1+2a,2+2a\}$,
$\{2+a,2+2a\}\}$.

\begin{displaymath}
\def\objectstyle{\scriptscriptstyle}
\xy /r5pc/:, {\xypolygon9"B"{\bullet}}, "B1"!{+U*++!L{\textrm{\normalsize0}}}
,"B2"!{+RD*+!LD{\textrm{\normalsize1}}},"B3"!{+RD*++!D{\textrm{\normalsize2}}}, 
"B4"!{+LDD*+!DR{\textrm{\normalsize2+{\it a}}}},"B5"!{+LD*++!R{\textrm{\normalsize {\it a}}}},
"B6"!{+LDD*+!UR{\textrm{\normalsize1+{\it a}}}},"B7"!{+LD*++!UR{\textrm{\normalsize1+2{\it a}}}},
"B8"!{+LDD*++!UL{\textrm{\normalsize2+2{\it a}}}},"B9"!{+LD*+!UL{\textrm{\normalsize2{\it a}}}},
"B3";"B1"**@{-},"B4";"B6"**@{-},"B7";"B9"**@{-},"B2";"B6"**@{-},"B2";"B7"**@{-},
"B5";"B1"**@{-},"B5";"B9"**@{-},"B8";"B3"**@{-},"B8";"B4"**@{-}
\endxy
\end{displaymath}
\begin{fig}
\begin{center}The Paley graph of order $9$
\end{center}
\end{fig}
\noindent{\it Example 3:}
Let $P=(V, E)$ be the Paley graph of order $13$. Here we have $p=13$, $n=1$ then $V(P)=\Z_{13}=\{0,1,2,3,4,5,6,7,8,9,10,11,12\}$ and $(\Z^*_{13})^2=\{1,3,4,9,10,12\}$. It follows that each vertex $x$ in $V(P)$ is adjacent exactly to $6$ vertices $x+1$, $x+3$, $x+4$, $x+9$, $x+10$, and $x+12$. So $E(P)=\{\{x,x+1\},\{x,x+3\},\{x,x+4\},\{x,x+9\},\{x,x+10\},\{x,x+12\}\ \forall x\in \Z_{13}\}$.

\begin{displaymath}
\def\objectstyle{\scriptscriptstyle}
\xy /r7pc/:, {\xypolygon13"B"{\bullet}}, "B1"!{+U*++!L{\textrm{\normalsize0}}}
,"B2"!{+RD*++!LD{\textrm{\normalsize1}}},"B3"!{+RD*++!LD{\textrm{\normalsize2}}}, 
"B4"!{+LDD*++!D{\textrm{\normalsize3}}},"B5"!{+LD*++!DR{\textrm{\normalsize 4}}},
"B6"!{+LDD*++!DR{\textrm{\normalsize5}}},"B7"!{+LD*++!R{\textrm{\normalsize6}}},
"B8"!{+LDD*++!R{\textrm{\normalsize7}}},"B9"!{+LD*++!UR{\textrm{\normalsize8}}},
"B10"!{+LDD*++!UR{\textrm{\normalsize9}}},"B11"!{+LD*++!U{\textrm{\normalsize10}}},
"B12"!{+LDD*++!UL{\textrm{\normalsize11}}},"B13"!{+LD*++!UL{\textrm{\normalsize12}}},
"B4";"B1"**@{-},"B4";"B7"**@{-},"B7";"B10"**@{-},"B10";"B13"**@{-},"B3";"B6"**@{-},
"B6";"B9"**@{-},"B9";"B12"**@{-},"B12";"B2"**@{-},"B2";"B5"**@{-},"B5";"B8"**@{-}
,"B8";"B11"**@{-},"B11";"B1"**@{-},
"B1";"B5"**@{-},"B9";"B5"**@{-},"B9";"B13"**@{-},"B13";"B4"**@{-},"B4";"B8"**@{-}
,"B8";"B12"**@{-},"B12";"B3"**@{-},"B3";"B7"**@{-},"B7";"B11"**@{-},"B11";"B2"**@{-}
,"B2";"B6"**@{-},"B6";"B10"**@{-},"B10";"B1"**@{-}
\endxy
\end{displaymath}
\begin{fig}
\begin{center}The Paley graph of order $13$
\end{center}
\end{fig}
\end{subsection}
\begin{subsection}{Properties}
In the previous examples we can see that the Paley graphs of order 5, 9, and 13 are connected, symmetric, self-complementary, and strongly regular. 

The following Propositions prove that these properties are true for every order $q$.
\begin{prop}
\label{Paley Symm}
The Paley graphs are symmetric.
\end{prop}
\proof
Let $P$ be the Paley graph of order $q = p^n$. To prove that $P$ is symmetric we need to prove that the automorphism group $Aut(P)$ acts transitively on $V(P)$ and $E(P)$. In other words, we need to prove that 

for every two vertices $x,y \in V(P)$ there exists $\phi \in Aut(P)$ such that $\phi(x)=y$, and 

for every two edges $\{x_1,y_1\},\{x_2,y_2\} \in E(P)$ there exists $\theta \in Aut(P)$ such that $\theta(x_1)=x_2,\theta(y_1)=y_2$. 
\\[10pt]
\indent Fix $a,b \in V(P)$ with $a \in (\F^{*}_{p^n})^2$ and define the nontrivial function $$\phi : V(P) \rightarrow V(P) {\textrm \ with \ } \phi(x)=ax+b \ \forall x\in V(P).$$ 
\indent We show that $\phi$ is an automorphism. Easily, we can see that $\phi$ is one to one, because $$\phi(x_1)-\phi(x_2)=0\Leftrightarrow (ax_1+b)-(ax_2+b)=0\Leftrightarrow a(x_1-x_2)+b-b=0$$ 
$$\Leftrightarrow x_1-x_2=0.$$ 
Since for every $y \in V(P)$, we have $$a^{-1}y-a^{-1}b=x \in V(P) \textrm{\ with \ }\phi(x)=a(a^{-1}y-a^{-1}b)+b=y.$$ 
Thus $\phi$ is onto. 
\\[10pt]
\indent Since $\{x,y\} \in E(P) \Leftrightarrow x-y \in (\F^{*}_{p^n})^2 \Leftrightarrow a(x-y)+b-b \in (\F^{*}_{p^n})^2 \Leftrightarrow (ax+b)-(ay+b) \in (\F^{*}_{p^n})^2 \Leftrightarrow \phi(x)-\phi(y) \in (\F^{*}_{p^n})^2 \Leftrightarrow \{\phi(x),\phi(y)\} \in E(P)$, this proves that $\phi \in Aut(P)$. 
\\[10pt]
\indent Moreover, for every two vertices $x,y \in V(P)$, take $a=1 \in (\F^{*}_{p^n})^2$ and $b=y-x \in V(P)$, the mapping $\phi : V(P) \rightarrow V(P)$ defined by $\phi(x)=ax+b$ is an automorphism with $\phi(x)=y$. Thus $Aut(P)$ acts transitively on $V(P)$. 
\\[10pt]
\indent Finally, for every two edges $\{x_1,y_1\},\{x_2,y_2\} \in E(P)$ we can find $$a=(x_2-y_2)(x_1-y_1)^{-1} \in (\F^{*}_{p^n})^2 {\textrm \ and \ } b=x_2-ax_1 \in V(P)$$ so that $\theta : V(P) \rightarrow V(P)$ with $\theta(x)=ax+b$ is an automorphism with $\theta(x_1)=x_2,\theta(y_1)=y_2$. Thus $Aut(P)$ acts transitively on $E(P)$.
\hspace{\stretch{1}} $\Box$
\vspace{10pt}
\begin{prop}
Let $P$ be the Paley graph of order $q = p^n$, then $P$ is a self-complementary graph.
\end{prop}
\proof 
Let $r$ be a quadratic nonresidue modulo $q$, consider the function $$f: V(P) \longrightarrow V(\overline{P})\textrm{ defined by }f(x)=rx.$$ 

The function $f$ is well defined, because 
$$\{x,y\} \in E(P) \Leftrightarrow (x-y) \in (\F^{*}_{p^n})^2 \Leftrightarrow f(x)-f(y)=rx-ry=r(x-y) \notin (\F^{*}_{p^n})^2$$ $$\Leftrightarrow \{f(x), f(y)\} \in E(\overline{P}).$$ 

Now we prove that $f$ is a bijection. Clearly, $f$ is injective, since $$(x-y)=0 \Leftrightarrow 0=r(x-y)=rx-ry=f(x)-f(y).$$  
\indent Since gcd $(r,q)=1$, there exist $$a , b \in \Z \textrm{ with }1=qa+rb \Leftrightarrow rb\equiv 1 \ (\textrm {mod } q).$$ 

Thus $f(bx)= rbx = x$, so $f$ is surjective. 
\hspace{\stretch{1}} $\Box$
\vspace{10pt}
\begin{prop}
 Let $P$ be the Paley graph of order $q = p^n$, then $P$ is a strongly regular graph with parameters 
\begin{displaymath}(q, \frac{q-1}{2}, \frac{q-5}{4}, \frac{q-1}{4}). \end{displaymath}
\end{prop}
\proof First, we prove that each vertex has degree $\frac{q-1}{2}$. 
\\[10pt]
\indent Fix $x \in V(P)$, we have $N(x)=\{z\in V(P) : x-z = s \in (\F^{*}_{p^n})^2\}$. If $x-z_1 = s, x-z_2 = s$ then $z_1= x-s = z_2$, so for all $s\in (\F^{*}_{p^n})^2$ there exists a unique $z\in V(P)$ such that $x-z = s$. 
\\[10pt]
\indent Thus there exists a one to one correspondence between the number of elements of $N(x)$ and the number of elements of $(\F^{*}_{p^n})^2$, so all vertices have the same degree $d(x)=|N(x)|=|(\F^{*}_{p^n})^2|$. 
\\[10pt]
\indent Now we calculate $|(\F^{*}_{p^n})^2|$, we have $|\F^{*}_{p^n}|=q-1$ and if $x\neq y \in \F^{*}_{p^n}$ then $x^2=y^2 \Leftrightarrow 0=x^2-y^2=(x-y)(x+y)\Leftrightarrow x=-y$. Thus $|(\F^{*}_{p^n})^2|=\frac{q-1}{2}$.
\\[10pt]
\indent Second, we prove that every two adjacent vertices have $\frac{q-5}{4}$ common neighbors and every two non-adjacent vertices have $\frac{q-1}{4}$ common neighbors. 
\\[10pt]
\indent Let $x \in V(P)=(\F^{*}_{p^n}), A= N(x), B= V(P)\setminus(A\cup\{x\})$. If $y \in A$ and $z \in B$, we want to prove that $|A\cap N(y)|=\frac{q-5}{4}$ and $|A\cap N(z)|=\frac{q-1}{4}$. 
\\[10pt]
\indent Because $P$ is symmetric, we can assume that there is an integer $l$ with every vertex $y \in A$ is joined to $l$ vertices in $B$ ($|N(y)\cap B|=l$). Moreover, because $P$ is self-complementary, every vertex $z \in B$ is not joined to $l$ vertices in $A$ ($|(V(P)\setminus N(z))\cap A|=l$). 
\\[10pt]
\indent To find $l$ we calculate $|A||B|$ from two sides. First $|A|= |N(x)|=\frac{q-1}{2}$, $|B|=|V(P)|-|A\cup\{x\}|=q-(\frac{q-1}{2}+1)=q-\frac{q+1}{2}=\frac{2q-q-1}{2}=\frac{q-1}{2}$, which means that $|A||B|=(\frac{q-1}{2})^2$. 
\\[10pt]
\indent Second $|A||B|=|A\times B|=|\{(a,b): a\in A,b\in B\textrm { and }\{a,b\} \in E(P)\}|+$ 
$|\{(a,b): a\in A,b\in B\textrm { and }\{a,b\} \notin E(P)\}|=\frac{q-1}{2}l+l\frac{q-1}{2}=2l\frac{q-1}{2}$. 
\\[10pt]
\indent So $|A||B|=(\frac{q-1}{2})^2=2l\frac{q-1}{2}$, which gives $2l=\frac{q-1}{2}$ or $l=\frac{q-1}{4}$. 
\\[10pt]
\indent Now we can calculate $|A\cap N(y)|$ and $|A\cap N(z)|$. Since 
\\[10pt]
$\frac{q-1}{2}=|N(y)|=|A\cap N(y)|+|B\cap N(y)|+|\{x\}\cap N(y)|=|A\cap N(y)|+l+1=
|A\cap N(y)|+\frac{q-1}{4}+1=|A\cap N(y)|+\frac{q+3}{4},$
\\[10pt] 
we have $|A\cap N(y)|=\frac{q-1}{2}-\frac{q+3}{4}=\frac{2q-2-q-3}{4}=\frac{q-5}{4}$.
\\[10pt]
\indent Since $\frac{q-1}{2}=|A|=|(V(P)\setminus N(z))\cap A|+|N(z)\cap A|=l+|N(z)\cap A|$ $=\frac{q-1}{4}+|N(z)\cap A|,$
\\[10pt] 
hence $|N(z)\cap A|=\frac{q-1}{2}-\frac{q-1}{4}=\frac{q-1}{4}.$
\\[10pt]
Then $P$ is a strongly regular graph with parameters 
$$(q, \frac{q-1}{2}, \frac{q-5}{4}, \frac{q-1}{4}).$$\hspace{\stretch{1}} $\Box$

\begin{cor}
The Paley graphs are connected.
\end{cor}
\proof
Let $P$ be the Paley graph of order $q = p^n$. Let $x,y$ be two vertices in $V(P)$, then $x,y$ are adjacent or non-adjacent. If $x,y$ are adjacent, then there exists a path of length $1$ connected $x$ and $y$. 
\\[10pt]
\indent If $x,y$ are non-adjacent, then $x$ and $y$ have at least one common neighbor $z$ because $q\geq 5$ means that $\frac{q-1}{4} \geq 1$. So there exists a path $xe_1ze_2y$ with $e_1=\{x,z\},e_2=\{z,y\}$ of length $2$ connected $x$ and $y$.
\\[10pt]
\indent Thus in all cases every two vertices in $V(P)$ are connected by a path.\hspace{\stretch{1}} $\Box$
\\[10pt]
\indent Now we know that the Paley graphs are self-complementary symmetric graphs. So the question now is: Are there any self-complementary symmetric graphs other than Paley graphs? 
\\[10pt]
\indent Peisert proved in \cite{W.P1} that the Paley graphs of prime order are the only self-complementary symmetric graphs of prime order and he proved in \cite{W.P2} that a graph $G$ is self-complementary and symmetric if and only if $ |G| = p^n$ for some prime $p,\ p^n \equiv 1$ (mod $4$), and $G$ is a Paley graph or a $\mathcal{P}^*$-graph or is isomorphic to the exceptional graph $G(23^2)$. 
\\[10pt]
\indent The $\mathcal{P}^*$-graph is a graph with $V(\mathcal{P}^*)=\F_{p^n}$ and two vertices are adjacent if their difference belongs to the set $M =\{ g^j : j \equiv 0 ,1$ (mod $4$)$\}$ , where $g$ is a primitive root of the field. The graph $G(23^2)$ has $23^2$ vertices and is described in Section 3 in \cite{W.P2}.
\\[10pt]
\indent One of the properties of the Paley graphs is the 3-existentially closed property. As in \cite{Hadamard M.}, for a fixed integer $n \geq 1$, a graph $G$ is $n$-existentially closed, if for every $n$-element subset $S$ of the vertices, and for every subset $T$ of $S$, there is a vertex $x \notin S$ which is joined to every vertex in $T$ and to no vertex in $S \setminus T$. 
\\[10pt]
\indent The $n$-existentially closed graphs were first studied in \cite{Erdos}, where they were called graphs with property $P (n)$. Ananchuen and Caccetta, in \cite{ontheAdjpaley}, proved that all Paley graphs with at least 29 vertices are 3-existentially closed, and before \cite{Hadamard M.} they were the only known examples of strongly regular 3-existentially closed graphs. Now in \cite{Hadamard M.} we can find a new infinite family of 3-existentially closed graphs, that are strongly regular but not Paley graphs. For further background on $n$-existentially closed graphs the reader is directed to \cite{on the Adjacency}.
\\[10pt]
\indent Another property of the Paley graphs is also interesting. To understand it we need the following definition.
\\[10pt]
\noindent {\bf Definition:} If $m,\ n \in \N \cup \{0\}$ and $k \in \N$, a graph $G$ is said to have the property $P(m,n,k)$, if for any disjoint subsets $A$ and $B$ of $V(G)$ with $|A|=m$ and $|B|=n$ there exist at least $k$ other vertices, each of which is adjacent to every vertex in $A$ but not adjacent to any vertex in $B$. 

The set of graphs which have the property $P(m,n,k)$ is denoted by $\mathcal{G}(m,n,k)$.
\\[10pt]  
In \cite{ontheAdjGpaley} and \cite{ongraphs}, it has been proved that the Paley graph 
$$P_q \in \mathcal{G}(1,n,k) \textrm{ for every } q >{\big (}(n-2)2^n+2{\big )}\sqrt{q}+(n+2k-1)2^n-2n-1;$$
$$P_q \in \mathcal{G}(n,n,k) \textrm{ for every } q >{\big (}(2n-3)2^{2n-1}+2{\big )}\sqrt{q}+(n+2k-1)2^{2n-1}-2n^2-1;$$
$$\textrm{and }P_q \in \mathcal{G}(m,n,k) \textrm{ for every } q >{\big (}(t-3)2^{t-1}+2{\big )}\sqrt{q}+(t+2k-1)2^{t-1}-1,$$
where $t \geq m+n.$ 
\end{subsection}
\end{section}
\chapter{Generalizations of The Paley Graphs}
\label{Chap3}
There are many generalizations of the Paley graphs. We will see some examples of these generalizations 
and some of its properties in the following section. In the second section we will define a new generalization, and we will study some of its properties in the third section.        
\begin{section}{Examples of some generalizations}
Since two vertices in the Paley graphs are adjacent if and only if their difference is a quadratic residue, we can generate other classes of graphs by using higher order residues. For example, in \cite{ontheAdjGpaley}, by using the cubic and quadruple residues Ananchuen has defined the cubic Paley graphs and the quadruple Paley graphs.
\begin{subsection}{The cubic and the quadruple Paley graphs}
\noindent {\bf Definition:} Let $q=p^n$ with odd prime $p$, $n \in \N$, and $q \equiv 1\ (\textrm{mod}\ 3)$. The graph $G^{(3)}_q$ with $$V(G^{(3)}_q)=\F_q \textrm{ and } E(G^{(3)}_q)=\{\{x,y\} : x,y \in \F_q,  x-y \in (\F^*_q)^3\}$$ is called {\it the cubic Paley graph}.
\\[10pt]
\indent Note that the set $E(G^{(3)}_q)$ in the definition is well defined because: $-1 = {-1}^3 \in (\F^*_q)^3$ implies that, 
$\{x,y\}$ is defined to be an edge if and only if $\{y,x\}$ is defined to be an edge.
\\[10pt]
\noindent {\bf Definition:} Let $q=p^n$ with odd prime $p$, $n \in \N$, and $q \equiv 1\ (\textrm{mod}\ 8)$. The graph $G^{(4)}_q$ with $$V(G^{(4)}_q)=\F_q \textrm{ and } E(G^{(4)}_q)=\{\{x,y\} : x,y \in \F_q,  x-y \in (\F^*_q)^4\}$$ is called {\it the quadruple Paley graph}.
\\[10pt]
\indent Note that the set $E(G^{(4)}_q)$ in the definition is well defined because : 
We have $q \equiv 1\ (\textrm{mod}\ 8)$, so $8 \mid (q -1)$. If $g$ is  a generator of the group $\F^*_{q}$ then $g^{\frac{q-1}{2}}=(g^{\frac{q-1}{8}})^4 =-1$ which means that $-1 \in ({\F_q}^*)^4$. Thus 
$\{x,y\}$ is defined to be an edge if and only if $\{y,x\}$ is defined to be an edge. 
\\[10pt]
\indent The following figure gives an example:
\begin{displaymath}
\def\objectstyle{\scriptscriptstyle}
\xy /r5pc/:, {\xypolygon13"B"{\bullet}}, "B1"!{+U*++!L{\textrm{\normalsize0}}}
,"B2"!{+RD*++!LD{\textrm{\normalsize1}}},"B3"!{+RD*++!LD{\textrm{\normalsize2}}}, 
"B4"!{+LDD*++!D{\textrm{\normalsize3}}},"B5"!{+LD*++!DR{\textrm{\normalsize 4}}},
"B6"!{+LDD*++!DR{\textrm{\normalsize5}}},"B7"!{+LD*++!R{\textrm{\normalsize6}}},
"B8"!{+LDD*++!R{\textrm{\normalsize7}}},"B9"!{+LD*++!UR{\textrm{\normalsize8}}},
"B10"!{+LDD*++!UR{\textrm{\normalsize9}}},"B11"!{+LD*++!U{\textrm{\normalsize10}}},
"B12"!{+LDD*++!UL{\textrm{\normalsize11}}},"B13"!{+LD*++!UL{\textrm{\normalsize12}}},
"B1";"B6"**@{-},"B6";"B11"**@{-},"B11";"B3"**@{-},"B3";"B8"**@{-},"B8";"B13"**@{-},
"B13";"B5"**@{-},"B5";"B10"**@{-},"B10";"B2"**@{-},"B2";"B7"**@{-},"B7";"B12"**@{-}
,"B12";"B4"**@{-},"B4";"B9"**@{-},
"B9";"B1"**@{-}
\endxy
\textrm{ \quad\quad\quad }
\xy /r5pc/:, {\xypolygon17"B"{\bullet}}, "B1"!{+U*++!L{\textrm{\normalsize0}}}
,"B2"!{+RD*++!LD{\textrm{\normalsize1}}},"B3"!{+RD*+!LD{\textrm{\normalsize2}}}, 
"B4"!{+LDD*++!D{\textrm{\normalsize3}}},"B5"!{+LD*++!D{\textrm{\normalsize4}}},
"B6"!{+LDD*++!D{\textrm{\normalsize5}}},"B7"!{+LD*+!DR{\textrm{\normalsize6}}},
"B8"!{+LDD*+!DR{\textrm{\normalsize7}}},"B9"!{+LD*++!R{\textrm{\normalsize8}}},
"B10"!{+LDD*++!R{\textrm{\normalsize9}}},"B11"!{+LD*++!R{\textrm{\normalsize10}}},
"B12"!{+LDD*+!UR{\textrm{\normalsize11}}},"B13"!{+LD*++!U{\textrm{\normalsize12}}},
"B14"!{+LD*++!U{\textrm{\normalsize13}}},"B15"!{+LD*+!UL{\textrm{\normalsize14}}},
"B16"!{+LD*+!UL{\textrm{\normalsize15}}},"B17"!{+LD*+!UL{\textrm{\normalsize16}}},
"B1";"B5"**@{-},"B5";"B9"**@{-},"B9";"B13"**@{-},"B13";"B17"**@{-},"B17";"B4"**@{-},
"B4";"B8"**@{-},"B8";"B12"**@{-},"B12";"B16"**@{-},"B16";"B3"**@{-},"B3";"B7"**@{-}
,"B7";"B11"**@{-},"B11";"B15"**@{-},"B15";"B2"**@{-},"B2";"B6"**@{-},"B6";"B10"**@{-},
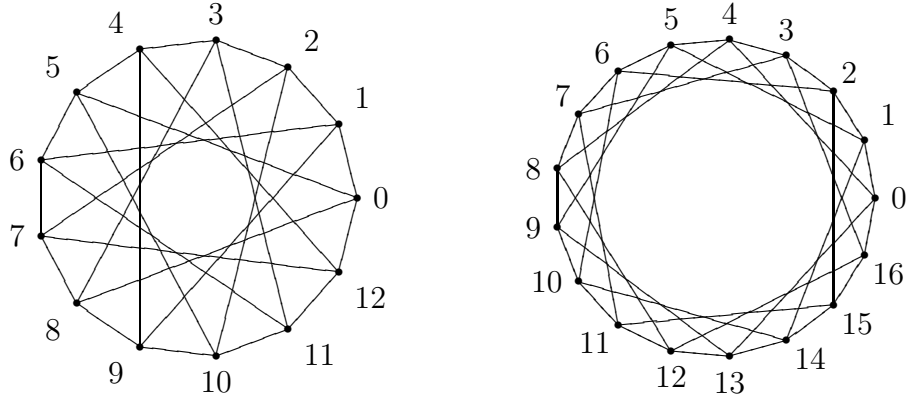
"B10";"B14"**@{-},"B14";"B1"**@{-}
\endxy
\end{displaymath}
\begin{fig}
\begin{center}The cubic Paley graph $G^{(3)}_{13}$ and the quadruple Paley graph $G^{(4)}_{17}$
\end{center}
\end{fig}
\vspace{10pt}
\noindent In \cite{ontheAdjGpaley}, Ananchuen has proved that the cubic Paley graphs 
$$G^{(3)}_q \in \mathcal{G}(2,2,k) \textrm{ for every } q > [\frac{1}{4}(79+3\sqrt{36k+701})]^2;$$ 
$$G^{(3)}_q \in \mathcal{G}(m,n,k) \textrm{ for every } q > (t2^{t-1}-2^t+1)2^m\sqrt{q}+(m+2n+3k-3)2^{-n}3^{t-1},$$
where $t \geq m+n$; and the quadruple Paley graphs
$$G^{(4)}_q \in \mathcal{G}(m,n,k) \textrm{ for every } q > (t2^{t-1}-2^t+1)3^m\sqrt{q}+(m+3n+4k-4)3^{-n}4^{t-1},$$
where $t \geq m+n$.
\\[10pt]
\noindent Ananchuen and Caccetta have proved in \cite{CubQuadnec} that the cubic Paley graphs are $n$-existentially closed whenever
$q \geq n^2 2^{4n-2}$ and the quadruple Paley graphs are $n$-existentially closed whenever $q \geq 9n^2 6^{2n-2}$.
\end{subsection}
\begin{subsection}{The generalized Paley graphs}
In \cite{onGenPal}, Lim and Praeger have defined the following generalization of the Paley graphs.
\\[10pt]
\noindent {\bf Definition:} Let $\F_q$ be a finite field of order $q$, and let $k$ be a divisor of $q - 1$ such that $k \geq 2$, and if $q$ is odd then $\frac{q-1}{k}$ is even. Let $S$ be the subgroup of order $\frac{q-1}{k}$ of the multiplicative group $\F^*_q$. Then the generalized Paley graph GPaley$(q,\frac{q-1}{k})$ is the graph with vertex set $\F_q$ and edges all pairs $\{x, y\}$ such that $x - y \in S$.
\\[10pt]
\indent Note that in the definition they require $\frac{q-1}{k}$ to be even when $q$ is odd, and hence in
all cases $S = -S$, so that the adjacency relation is symmetric ($x - y \in S$ if and only if $y - x \in S$). 

Also we can see that if $q \equiv 1\mod 4$ and $k=2$, then GPaley$(q,\frac{q-1}{k})$ is the Paley graph $P_q$. 
\\[10pt]
\indent As an example consider $q=11$ then $q-1=10$, and since $\frac{q-1}{k}$ should be even and $k \geq 2$, we have only one choice $k=5$. Then $|S|=\frac{10}{5}=2$ and $S=\{1,10\}$, it follows that GPaley$(11,2)$ is the cycle $C_{11}$.
\\[10pt]
\indent Moreover, they have studied in \cite{onGenPal} the automorphism groups of this generalized Paley graphs, and in some cases, compute their full automorphism groups. Moreover they have determined precisely when these graphs are connected.
\end{subsection}

\end{section}

\begin{section}{Definition and examples}
Now we will give a new generalization of the Paley graphs.
\\[10pt]
\noindent {\bf Definition:} Let $q=p^n$ with odd prime $p$, $n \in \N$, and $m \geq 3$ be an odd integer. We will denote by $m\textrm{-}P_q$, the graph with $V(m\textrm{-}P_q)=\F_q$ and $E(m\textrm{-}P_q)=\{\{x,y\} : x,y \in \F_q,  x-y \in ({\F_q}^*)^m\}$. Such a graph will be called $m$-Paley graph. 
\\[10pt]
\indent Note that the set $E(m\textrm{-}P_q)$ in the definition is well defined because : $-1 = {-1}^m \in ({\F_q}^*)^m$ implies that 
\begin{center}
$x-y\in (\F^*_{p^n})^m$ if and only if $y-x\in (\F^*_{p^n})^m$.  
\end{center}

\indent The list of integers which can be considered as the order of $m$-Paley graph starts with $3$,$5$,$7$,$9$,$11$,$13$,$17$,$19$,$23$,$25$,$27$,$29$,$31$. 

In the following examples, we will show the $m$-Paley graphs explicitly for the first three cases. 
\\[10pt]
{\it Example 1:}
The $m$-Paley graph of order $3$, for every odd integer $m \geq 3$, is the cycle $C_3$. 

In order to see that, let $m\textrm{-}P_3=(V, E)$ be the $m$-Paley graph of order $3$ then $V(m\textrm{-}P_3)=\Z_3=\{0,1,2\}$ and $(\Z^*_3)^m=\{1,2\}$. 

It follows that $E(m\textrm{-}P_3)=\{\{0,1\},\{1,2\},\{2,0\}\}.$
\\[10pt]
{\it Example 2:}
The $m$-Paley graph of order $5$, for every odd integer $m \geq 3$, is the complete graph $K_5$. 

In order to see that, let $m\textrm{-}P_5=(V, E)$ be the $m$-Paley graph of order $5$ and $m=2k+1$ with positive integer $k$, then $V(m\textrm{-}P_5)=\Z_5=\{0,1,2,3,4\}$ and $1^m=1$, $4^m={-1}^m=-1$. To calculate $2^m$, we have two cases:

If $k$ is odd then $2^m=2^{2k+1}=(2^2)^k 2=(-1)^k 2=3$ and if $k$ is even then $2^m=2^{2k+1}=(-1)^k 2=2$. 
To calculate $3^m$, we have also two cases:
 
If $k$ is odd then $3^m=3^{2k+1}=(3^2)^k 3=(-1)^k 3=2$ and if $k$ is even then $3^m=3^{2k+1}=(-1)^k 3=3$. 

Which gives $(\Z^*_5)^m=\{1,2,3,4\}$ for every odd integer $m$, it follows that $E(m\textrm{-}P_5)=\{\{0,1\},\{0,2\},\{0,3\},\{0,4\},\{1,2\},\{1,3\},\{1,4\},\{2,3\},\{2,4\},$ \\
$\{3,4\}\}$.
\\[10pt]
{\it Example 3:}
Let $m\textrm{-}P_7=(V, E)$ be the $m$-Paley graph of order $7$, then $V(m\textrm{-}P_7)=\Z_7=\{0,1,2,3,4,5,6\}$.

Consider $m=3$, then $1^3=1$, $2^3=1$, $3^3=6$, $4^3=1$, $5^3=6$, $6^3=6$. It follows that 
$(\Z^*_7)^3=\{1,6\}$ and 

\noindent $E(3\textrm{-}P_7)=\{\{0,1\},\{1,2\},\{2,3\},\{3,4\},\{4,5\},\{5,6\},\{6,0\}\}$. 

Consider $m=5$, then $1^5=1$, $2^5=4$, $3^5=5$, $4^5=2$, $5^5=3$, $6^5=6$. It follows that

\noindent $(\Z^*_7)^5=\{1,2,3,4,5,6\}$ and $5\textrm{-}P_7$ is the complete graph $K_7$. 

Consider $m=7$, then $1^7=1$, $2^7=2$, $3^7=3$, $4^7=4$, $5^7=5$, $6^7=6$. It follows that 

\noindent $(\Z^*_7)^7=\{1,2,3,4,5,6\}$ and $7\textrm{-}P_7$ is the complete graph $K_7$. 

Consider $m=9$, then $1^9=1$, $2^9=1$, $3^9=6$, $4^9=1$, $5^9=6$, $6^9=6$. It follows that 
$(\Z^*_7)^9=\{1,6\}$ and 

\noindent $E(9\textrm{-}P_7)=\{\{0,1\},\{1,2\},\{2,3\},\{3,4\},\{4,5\},\{5,6\},\{6,0\}\}$.
\begin{displaymath}
\def\objectstyle{\scriptscriptstyle}
\xy /r4pc/:, {\xypolygon7"B"{\bullet}}, "B1"!{+U*+!L{\textrm{\normalsize0}}}
,"B2"!{+RD*+!D{\textrm{\normalsize1}}}
,"B3"!{+RD*+!R{\textrm{\normalsize2}}}, "B4"!{+LDD*+!R{\textrm{\normalsize3}}}
,"B5"!{+RD*+!RU{\textrm{\normalsize4}}}, "B6"!{+LDD*+!UL{\textrm{\normalsize5}}}
,"B7"!{+RD*+!L{\textrm{\normalsize6}}}
\endxy
\textrm{ \quad\quad\quad\quad\quad\quad }
\xy /r4pc/:, {\xypolygon7"B"{\bullet}}, "B1"!{+U*+!L{\textrm{\normalsize0}}}
,"B2"!{+RD*+!D{\textrm{\normalsize1}}}
,"B3"!{+RD*+!R{\textrm{\normalsize2}}}, "B4"!{+LDD*+!R{\textrm{\normalsize3}}}
,"B5"!{+RD*+!RU{\textrm{\normalsize4}}}, "B6"!{+LDD*+!UL{\textrm{\normalsize5}}}
,"B7"!{+RD*+!L{\textrm{\normalsize6}}},
"B1";"B2"**@{-},"B1";"B3"**@{-},"B1";"B4"**@{-},"B1";"B5"**@{-},"B1";"B6"**@{-},"B1";"B7"**@{-},
"B2";"B3"**@{-},"B2";"B4"**@{-},"B2";"B5"**@{-},"B2";"B6"**@{-},"B2";"B7"**@{-},
"B3";"B4"**@{-},"B3";"B5"**@{-},"B3";"B6"**@{-},"B3";"B7"**@{-},
"B4";"B5"**@{-},"B4";"B6"**@{-},"B4";"B7"**@{-},
"B5";"B6"**@{-},"B5";"B7"**@{-},
"B6";"B7"**@{-}
\endxy
\end{displaymath}
\begin{center}
$3\textrm{-}P_7$, $9\textrm{-}P_7$ \quad\quad\quad\quad\quad
\quad\quad\quad\quad\quad $5\textrm{-}P_7$, $7\textrm{-}P_7$\end{center}
\begin{fig}
\begin{center}
\end{center}
\end{fig}
\noindent We can see that if gcd $(m,6=7-1)= 1$ then $m\textrm{-}P_7$ is the complete graph $K_7$, and if gcd $(m,6=7-1)= 3$ then $m\textrm{-}P_7$ is the cycle $C_7$.
\end{section}
\begin{section}{Properties}
We have seen in the previous examples that the $m$-Paley graphs $m\textrm{-}P_q$ are the complete graph $K_q$ in some cases and in other cases are not. The question, which we consider now, is : When are the $m$-Paley graphs complete? The answer to this question can be found in the following theorem.
\begin{theo}
\label{mPaleycom}
Let $m\textrm{-}P_q=(V,E)$ be the $m$-Paley graph of order $q$ and $d=$ gcd $(m,q-1)$, then 
\begin{center}
$m\textrm{-}P_q$ is complete if and only if $d=1$. 
\end{center}
\end{theo}
\proof
$(\Rightarrow)$ If $m\textrm{-}P_q$ is complete then $d=1$: Assume that $m\textrm{-}P_q$ is complete, then for all $x\neq y \in \F_q$ we have $\{x,y\} \in E(m\textrm{-}P_q)$, it follows $x-y \in ({\F^*_q})^m$.

Clearly $(\F^*_q)^m \subseteq \F^*_q$. If $a \in \F^*_q $ then we can find $x,\ y \in \F_q$ with $a=x-y \in ({\F^*_q})^m$. It follows that $(\F^*_q)^m=\F^*_q$, which means that for all $a \in {\F^*_q}$ there exists $b \in {\F^*_q}$ with $a=b^m$. 
\\[10pt]
\indent By Theorem \ref{TH:cyclic} we know that ${\F^*_q}$ is cyclic. Let $g$ be a generator of ${\F^*_q}$, then for all $i \in \{1,2,\dots,q-1\}$ there exists $j \in \{1,2,\dots,q-1\}$ with $g^i=(g^j)^m$ or $g^{i-jm}=1$. Since order $g$ is $q-1$, we have $q-1\mid i-jm$, so $i-mj=(q-1)k$ for some $k$. Thus $i=mj+(q-1)k$.
\\[10pt]
\indent Now to prove that $d=1$, assume the contrary. Let $d > 1$, then $$m=dm_1, \ (q-1)=dk_1{\textrm \ for\ some\ } m_1,k_1.$$ 
\indent So we can write $i$ in the form $i=dm_1j+dk_1k=d(m_1j+k_1k)$, which means that $i$ must be a multiple of $d$, but $i$ is an arbitrary element of the set $\{1,2,\dots,q-1\}$. Thus we have a contradiction, then $d=1$.
\\[10pt]
\noindent$(\Leftarrow)$ If $d=1$ then $m\textrm{-}P_q$ is complete: Assume that $d=1$ and $m\textrm{-}P_q$ is not complete, then in the set $\{g^m,(g^2)^m,\dots,(g^{q-1})^m=1\}$ there exists two equal elements. Let $(g^i)^m=(g^j)^m$ with $i,j \in \{1,2,\dots,q-1\}, i>j$, then $g^{(i-j)m}=1$, so $q-1\mid(i-j)m$. 
\\[10pt]
\indent From the assumption $\gcd (m,q-1)=d=1 \Rightarrow q-1\mid i-j$, which is impossible because $i-j < q-1$. Then $m\textrm{-}P_q$ is complete.\hspace{\stretch{1}} $\Box$

\vspace{10pt}

\noindent Note that the completeness of the $m$-Paley graph means that $(\F^*_q)^m=\F^*_q$. So as an application of Theorem \ref{mPaleycom} we have the following corollary.
\vspace{10pt}
\begin{cor}
In the field $\F_q$, the equation $x^m=a$, $a \in \F_q$ has exactly one solution if and only if $d=$ gcd $(m,q-1)=1$.  
\end{cor}
\vspace{10pt}
\noindent Note that $d$ is an odd integer because $m$ is an odd integer. Now the question is: How does the $m$-Paley graphs, if gcd $(m,q-1)=d > 1$, look like? The following proposition has the answer of this question.
\vspace{10pt}
\begin{prop}
\label{mPaleyreg}
Let $m\textrm{-}P_q=(V,E)$ be the $m$-Paley graph of order $q$,
\begin{center}
if $d=$ gcd $(m,q-1) > 1$, then $m\textrm{-}P_q$ is $\frac{q-1}{d}$-regular. 
\end{center}
\end{prop}
\proof
Let $x$ be any vertex in $V(m\textrm{-}P_q)$, then $y \in V(m\textrm{-}P_q)$ is adjacent to $x$ if and only if there exists $z \in (\F^*_q)^m$ with $x-y=z$. Which means that $|N(x)|=|(\F^*_q)^m|$ for all $x \in V(m\textrm{-}P_q)$, then $m\textrm{-}P_q$ is $|(\F^*_q)^m|$-regular. 
\\[10pt]
\indent Let us now prove that $|(\F^*_q)^m|=\frac{q-1}{d}$. Let $g$ be a generator of the group $\F^*_q$ and $(g^i)^m=(g^j)^m$ for some $i,j \in \{1,2,\dots,q-1\}$, then $g^{(i-j)m}=1$, which implies that $q-1\mid(i-j)m$. 
\\[10pt]
\indent Since $d\mid q-1,\  d\mid m$, it follows that $\frac{q-1}{d} \mid \frac{m}{d}(i-j)$. Since gcd $(\frac{q-1}{d},\frac{m}{d})=1$, we have $\frac{q-1}{d} \mid (i-j) \Leftrightarrow i\equiv j$ (mod $\frac{q-1}{d}$). So the following $d$ elements $$g^i,g^{i+\frac{q-1}{d}},g^{i+2\frac{q-1}{d}},\dots,g^{i+(d-1)\frac{q-1}{d}}$$ are the same. Thus $|(\F^*_q)^m|=\frac{q-1}{d}$ and $m\textrm{-}P_q$ is $\frac{q-1}{d}$-regular.\hspace{\stretch{1}} $\Box$
\vspace{10pt}
\begin{cor}
In the field $\F_q$, if $d=$ gcd $(m,q-1) > 1$, then the equation $x^m=a$, $a \in \F^*_q$ has exactly $d$ solutions.  
\end{cor}
\noindent As a special case, if $d=m$ and $a=1$ we get Lagrange's lemma.

\vspace{10pt}
We can see that the $m$-Paley graphs are not strongly regular in general. In example $3$ we have seen that $3\textrm{-}P_7$, $9\textrm{-}P_7$ are $C_7$ which is not strongly regular. So the question now is: Are the $m$-Paley graphs symmetric or self-complementary? 
\vspace{10pt}
\begin{prop}
The $m$-Paley graphs are symmetric.
\end{prop}
\proof
By the same proof of Proposition \ref{Paley Symm} with replacing $(\F^*_q)^2$ by $(\F^*_q)^m$ we get the result.\hspace{\stretch{1}} $\Box$
\vspace{10pt}
\begin{prop}
The $m$-Paley graphs are not self-complementary.
\end{prop}
\proof
Clearly, a self-complementary graph of order $q$ should have $\frac{q(q-1)}{4}$ edges. Let $m\textrm{-}P_q=(V,E)$ be the $m$-Paley graph of order $q$. 

If $d=$gcd$(m,q-1) > 1$, then by Proposition \ref{mPaleyreg} 
\begin{displaymath}
|E(m\textrm{-}P_q)|=\frac{1}{2} \sum_{x\in V(m\textrm{-}P_q)} d(x) =\frac{1}{2}q\frac{q-1}{d}. 
\end{displaymath}
Since $d\geqslant 3$, it follows that $|E(m\textrm{-}P_q)|<\frac{q(q-1)}{4}$. 

If $d=$ gcd $(m,q-1)=1$, then by Proposition \ref{mPaleycom} 
$$|E(m\textrm{-}P_q)|=q\frac{q-1}{2}>\frac{q(q-1)}{4}.$$ 

Thus the $m$-Paley graphs are not self-complementary.\hspace{\stretch{1}} $\Box$

\vspace{10pt}

\noindent Now let us ask the following question: Are the $m$-Paley graphs connected? 

\noindent Clearly, the $m$-Paley graph of order $q$ with $d= \gcd (m,q-1) = 1$, which is the complete graph $K_q$, is connected. So the case which we will study is the $m$-Paley graph of order $q$ with $d=\gcd (m,q-1) > 1$. 
\\[10pt]
\noindent Note that $d \neq q-1$ because $q-1$ is even and $d$ must be odd. So 

\hspace{100pt} $\frac{q-1}{2} \geq d=$gcd$(m,q-1) \geq 1$.
\vspace{10pt}
\begin{prop}
\label{Prop:prim}
Let $m\textrm{-}P_q=(V,E)$ be the $m$-Paley graph of order $q$,
if $d=$ \normalfont{gcd} $(m,q-1) > 1$ and $q$ is prime, then $m\textrm{-}P_q$ is connected. 
\end{prop}
\proof
Since $q$ is prime, $\F_q = \Z_q$. For every $x < y \in \F_q$, we have the sequence $x,x+1,x+2,\dots,x+(y-x-1)=y-1,x+(y-x)=y \in \F_q$. 
So $\{x,x+1\},\{x+1,x+2\},\dots,\{y-1,y\} \in E(m\textrm{-}P_q)$, because $1\in (\F^*_q)^m$ for each odd integer $m$. Thus $x\{x,x+1\}x+1\{x+1,x+2\}x+2\cdots y-1\{y-1,y\}y$ is a path in $m\textrm{-}P_q$ between $x$ and $y$. So $m\textrm{-}P_q$ is connected. 
\hspace{\stretch{1}} $\Box$
\vspace{10pt}

\noindent Note that as a special case of Proposition \ref{Prop:prim}, if $d=\gcd(m,q-1)= \frac{q-1}{2}$ and $q$ is prime, then $m\textrm{-}P_q$ is the cycle $C_q$, because $m\textrm{-}P_q$ is connected and $2$-regular. 

\begin{prop}
Let $m\textrm{-}P_q=(V,E)$ be the $m$-Paley graph of order $q$,
if $d=\gcd(m,q-1)= \frac{q-1}{2}$ and $q$ is not prime, then $m\textrm{-}P_q$ is disconnected. 
\end{prop}
\proof
Let $q=p^n, n > 1$, then $\F_q = \Z_p[x] /(f(x))$ where $f(x)$ is an irreducible polynomial of degree $n$ over $\Z_p$.
\\[10pt]
\indent Since $1 \in (\F^*_q)^m$, the path $0\{0,1\}1\{1,2\}2\cdots p-1\{p-1,0\}0\ $ in $m\textrm{-}P_q$ form a cycle, say $C_p$. By using Proposition \ref{mPaleyreg} with $d=$ gcd $(m,q-1)= \frac{q-1}{2}$, it follows that $m\textrm{-}P_q$ is $2$-regular. 
\\[10pt]
\indent Since $n > 1$, it follows that the set $\F_q \setminus\Z_p$ is not empty. Let $a \in \F_q \setminus\Z_p$, then $a$ is not adjacent to any vertex in $C_p$. So we cannot find a path in $m\textrm{-}P_q$ between any vertex in $\F_q \setminus\Z_p$ and any vertex in $C_p$. 

Thus $m\textrm{-}P_q$ is disconnected.\hspace{\stretch{1}} $\Box$
\\[20pt]
{\it Example:} Take $q=27, \ m=13$, then $d=\gcd(13,27-1)=13=\frac{27-1}{2}$. The following figure shows that the $13\textrm{-}P_{27}$ is disconnected 
\begin{displaymath}
\def\objectstyle{\scriptscriptstyle}
\xy /r6pc/:, {\xypolygon27"B"{~>{}\bullet}}, "B1"!{+U*+!L{\scriptstyle 0}}
,"B2"!{+RD*+!L{\scriptstyle 1}}
,"B3"!{+RD*+!L{\scriptstyle 2}}, "B4"!{+LDD*+!L{\scriptstyle a}}
,"B5"!{+RD*+!LD{\scriptstyle a+1}}, "B6"!{+LDD*+!LD{\scriptstyle a+2}}
,"B7"!{+RD*++!LD{\scriptstyle a^2}},
"B8"!{+RD*++!D{\scriptstyle a^2+1}}, "B9"!{+LDD*+!RD{\scriptstyle a^2+2}},
"B10"!{+RD*++!R{\scriptstyle 2a}}, "B11"!{+LDD*+!R{\scriptstyle 2a+1}},
"B12"!{+RD*+!R{\scriptstyle 2a+2}}, "B13"!{+LDD*+!R{\scriptstyle 2a^2}},
"B14"!{+RD*+!R{\scriptstyle 2a^2+1}}, "B15"!{+LDD*+!R{\scriptstyle 2a^2+2}},
"B16"!{+RD*+!R{\scriptstyle a^2+a}}, "B17"!{+LDD*+!R{\scriptstyle a^2+a+1}},
"B18"!{+RD*+!R{\scriptstyle a^2+a+2}}, "B19"!{+LDD*+!R{\scriptstyle a^2+2a}},
"B20"!{+RD*!RU{\scriptstyle a^2+2a+1}}, "B21"!{+LDD*++!U{\scriptstyle a^2+2a+2}},
"B22"!{+RD*!LU{\scriptstyle 2a^2+a}}, "B23"!{+LDD*++!L{\scriptstyle 2a^2+a+1}},
"B24"!{+RD*+!L{\scriptstyle 2a^2+a+2}}, "B25"!{+LDD*+!L{\scriptstyle 2a^2+2a}},
"B26"!{+RD*+!L{\scriptstyle 2a^2+2a+1}}, "B27"!{+LDD*+!L{\scriptstyle 2a^2+2a+2}},
"B3";"B1"**@{-},"B3";"B2"**@{-},"B1";"B2"**@{-},"B5";"B4"**@{-},"B6";"B4"**@{-},"B5";"B6"**@{-},
"B7";"B8"**@{-},"B7";"B9"**@{-},"B8";"B9"**@{-},"B10";"B11"**@{-},"B10";"B12"**@{-},"B11";"B12"**@{-},
"B13";"B14"**@{-},"B13";"B15"**@{-},"B15";"B14"**@{-},"B16";"B17"**@{-},"B16";"B18"**@{-},"B17";"B18"**@{-},
"B19";"B20"**@{-},"B19";"B21"**@{-},"B20";"B21"**@{-},"B22";"B23"**@{-},"B23";"B24"**@{-},"B22";"B24"**@{-},
"B25";"B26"**@{-},"B25";"B27"**@{-},"B26";"B27"**@{-},
\endxy
\end{displaymath}
\begin{fig}
\begin{center}
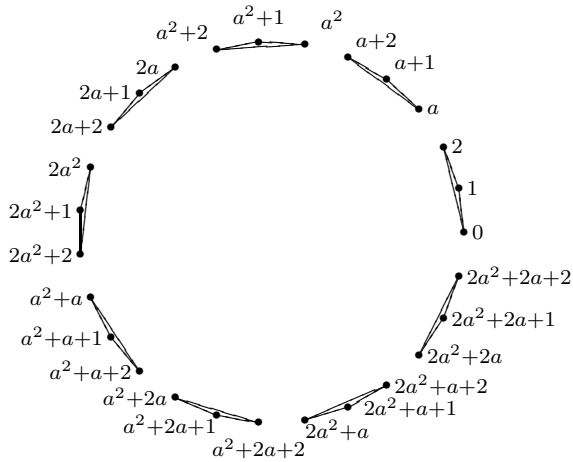
The $13\textrm{-}P_{27}$ graph, where $a$ is a root of an irreducible polynomial $f(x)$ of degree $3$ over $\Z_3$
\end{center}
\end{fig}

\noindent So now we know that the $m$-Paley graphs are not always connected and the case : $q= p^n,\ n > 1 $ and $ 1 < d=\gcd (m,q-1) < \frac{q-1}{2}$ is still open.

\end{section}

\newpage

\newpage
\thispagestyle{empty}


\noindent Hiermit versichere ich, die vorliegende Arbeit selbstst\"{a}ndig verfasst und keine anderen als die angegebenen Quellen und Hilfsmittel benutzt zu haben.
\vspace{2em}

\noindent D\"{u}sseldorf, den 27.5.2009\\

\vspace{2em}

\noindent\rule{0.3\textwidth}{0.4pt}\\
Ahmed Elsawy

\end{document}